\documentclass[12pt,leqno,fleqn,epsfig]{article}
\usepackage{amssymb,epsfig, amsmath, amsthm}
\usepackage{mathrsfs}       

\newcommand{\R}{\mathbb{R}}

\newcommand{\N}{\mathbb{N}}

\newcommand{\ba}{\boldsymbol a}
\newcommand{\bb}{\boldsymbol b}
\newcommand{\bc}{\boldsymbol c}
\newcommand{\bd}{\boldsymbol d}
\newcommand{\bfe}{\boldsymbol e}
\newcommand{\bbf}{\boldsymbol f}
\newcommand{\bg}{\boldsymbol g}
\newcommand{\bh}{\boldsymbol h}
\newcommand{\bi}{\boldsymbol i}
\newcommand{\bj}{\boldsymbol j}
\newcommand{\bk}{\boldsymbol k}
\newcommand{\bl}{\boldsymbol l}
\newcommand{\bm}{\boldsymbol m}
\newcommand{\bn}{\boldsymbol n}
\newcommand{\bo}{\boldsymbol o}
\newcommand{\bp}{\boldsymbol p}
\newcommand{\bq}{\boldsymbol q}
\newcommand{\br}{\boldsymbol r}
\newcommand{\bs}{\boldsymbol s}
\newcommand{\bt}{\boldsymbol t}
\newcommand{\bu}{\boldsymbol u}
\newcommand{\bv}{\boldsymbol v}
\newcommand{\bw}{\boldsymbol w} 
\newcommand{\bx}{\boldsymbol x}
\newcommand{\by}{\boldsymbol y}
\newcommand{\bz}{\boldsymbol z}
\newcommand{\buu}{ \underline{\boldsymbol u}}

\newcommand{\bA}{\boldsymbol A}
\newcommand{\bB}{\boldsymbol B}
\newcommand{\bC}{\boldsymbol C}
\newcommand{\bD}{\boldsymbol D}
\newcommand{\bE}{\boldsymbol E}
\newcommand{\bF}{\boldsymbol F}
\newcommand{\bG}{\boldsymbol G}
\newcommand{\bH}{\boldsymbol H}
\newcommand{\bI}{\boldsymbol I}
\newcommand{\bJ}{\boldsymbol J}
\newcommand{\bK}{\boldsymbol K}
\newcommand{\bL}{\boldsymbol L}
\newcommand{\bM}{\boldsymbol M}
\newcommand{\bO}{\boldsymbol O}
\newcommand{\bP}{\boldsymbol P}
\newcommand{\bQ}{\boldsymbol Q}
\newcommand{\bR}{\boldsymbol R}
\newcommand{\bS}{\boldsymbol S}
\newcommand{\bT}{\boldsymbol T}
\newcommand{\bU}{\boldsymbol U}
\newcommand{\bV}{\boldsymbol V}
\newcommand{\bW}{\boldsymbol W}
\newcommand{\bX}{\boldsymbol X}
\newcommand{\bY}{\boldsymbol Y}
\newcommand{\bZ}{\boldsymbol Z}

\newcommand{\bfvarphi}{\boldsymbol\varphi}

\newcommand{\bfpsi}{\boldsymbol\psi}

\newcommand{\bPhi}{\boldsymbol \Phi}
\newcommand{\bfomega}{\boldsymbol \omega}

\newcommand{\bLambda}{\boldsymbol \Lambda}

\newcommand{\mes}{\operatorname{\rm mes}}    
\newcommand{\esssup}{\operatorname*{ess\,sup}}

\newcommand{\const}{\operatorname*{const}}

\newcommand{\be}{\begin{equation}}
\newcommand{\ee}{\end{equation}}
\newcommand{\bea}{\begin{eqnarray}}
\newcommand{\eea}{\end{eqnarray}}
\newcommand{\bean}{\begin{eqnarray*}}
\newcommand{\eean}{\end{eqnarray*}}

\newcommand{\var}{\varepsilon}

\newcommand{\intl}{\int\limits}
\newcommand{\suml}{\sum\limits}

\newcommand{\Beweisende}{\rule{0.2cm}{0.2cm}}

\newcommand{\D}{\displaystyle}

\newcommand{\intmw}{{\int\hspace{-830000sp}-\!\!}}

\newcounter{secnum}

\newtheorem{thm}{Theorem}[section]

\newtheorem{lem}[thm]{Lemma}

\theoremstyle{definition}
\newtheorem{defin}[thm]{Definition}
\newtheorem{rem}[thm]{Remark}

\title{On partial regularity for the $3D$ non-stationary  Hall magnetohydrodynamics equations on the plane}
 \author{Dongho Chae$^*$  and J\"{o}rg Wolf $^\dagger$\\
\ \\
 $*$Department of Mathematics\\
Chung-Ang University\\
 Seoul 156-756, Republic of Korea\\
 e-mail: dchae@cau.ac.kr\\
and \\
$\dagger$Department of Financial Engineering\\
 Ajou University\\                                           
Suwon 443-749, Republic of Korea\\
e-mail: jwolf@math.hu-berlin.de}
\date{}
\begin{document}
\maketitle
\begin{abstract}
We study partial regularity of  weak solutions of the  3D valued non-stationary Hall magnetohydrodynamics equations on $ \Bbb R^2$.  In particular we prove the existence 
of a weak solution whose set of possible singularities has the space-time Hausdorff dimension at most two. 
  \noindent{\bf
keywords:}  non-stationary Hall-MHD equations, partial regularity

\end{abstract}



\section{Introduction and the main theorem}
\label{sec:-1}
\setcounter{secnum}{\value{section} \setcounter{equation}{0}
\renewcommand{\theequation}{\mbox{\arabic{secnum}.\arabic{equation}}}}

We consider the homogeneous incompressible Hall magnetohydrodynamics(Hall-MHD) equations:
$$
\left\{ \aligned 
&\frac{\partial\bu}{\partial t} +\bu\cdot\nabla \bu+\nabla p=(\nabla\times \bB)\times \bB +\nu\Delta \bu +\bbf,\\
&\frac{\partial\bB}{\partial t}-\nabla\times (\bu\times \bB)+\nabla\times ((\nabla\times \bB)\times \bB)
=\mu \Delta B +\nabla \times \bg,\\
&\nabla \cdot \bu=0,\,\,\nabla \cdot \bB=0,\\ \endaligned \right.
$$
where the three dimensional vector fields $\bu=\bu(x,t)$ and $\bB=\bB(x,t)$ are the fluid velocity and the magnetic field respectively.  The scalar field $p=p(x,t)$ is the pressure, while the positive constants $\nu$ and $\mu$ represent the viscosity and the magnetic resistivity respectively. The given vector fields $\bbf$ and $\nabla \times \bg$ are  external forces  on the magnetically charged fluid flows.  
The  system has been studied first by Lighthill \cite{lig} in 1960.  We notice that comparing with the usual MHD system,  the  Hall-MHD system contains the extra term $\nabla\times ((\nabla\times \bB)\times \bB)$, called the Hall term.  The inclusion of this  term is essential in understanding the phenomena of magnetic reconnection, meaning the change of the topology of magnetic field lines. 
This  is observed in real physical situations such as  space plasma \cite{for, hom}, star formation \cite{war} and neutron star \cite{sha}.  
There are also many other physical phenomena that  requires the Hall-MHD system to describe  them (see e.g. \cite{miu, sim, pol} and the references therein).  
The Hall term is quadratically nonlinear, containing the second order derivative, and it causes major difficulties in the mathematical study of the Hall-MHD system.  Thanks to the orthogonality in $L^2$ of the Hall term with $\bB$, however, the energy inequality similar to the usual MHD case holds true. Using this fact the construction of the global in time weak solution can be achieved  without any difficulties, as has been shown in \cite{ach}.
Observing similar cancellation properties of the Hall term, the local in time well-posedness as well as the global in time well-posedness for small initial data was also established in \cite{cha1}, and has been refined in \cite{cha2}. 
Regarding the  question of energy conservation for weak solutions in the inviscid case we refer to \cite{dum}. 
For a special form of axially symmetric initial data the authors of \cite{fan} proved the global in time existence of classical solutions to the system. On the other hand,  the optimal temporal decay estimates are obtained  in \cite{cha3}. \\

\hspace{0.5cm}
Concerning the regularity of weak solutions, one can expect that the problem is more difficult than the Navier-Stokes equations and the usual MHD system.  Even the problem of regularity of stationary weak solution  has essential difficulty with current methods of the regularity theory, which is contrary to the case of the stationary Navier-Stokes equations. 
The partial regularity of stationary weak solutions has been obtained  recently by the current authors (cf. \cite{cha4}).  In the present  paper we investigate the partial regularity of weak solutions of {\em the non-stationary system}.  For the Navier-Stokes equations there are  many publications on this direction of study (see e.g. \cite{sch, caf, lad, lin, wol1}).  In the case of the 3D Hall-MHD system in $\Bbb R^3$, however, we encounter again essential difficulties  in constructing suitable weak solutions, satisfying desired form of localized energy inequality.

\hspace{0.5cm}\
In the current paper we focus on the case of 3D valued Hall-MHD system on the plane, which  is sometimes called the $2\frac12$ dimensional system. Physically the situation corresponds to the full 3D system having the translational symmetry in the $ x_3$ direction.
In this case, as will be shown in detail below,  although we cannot construct  suitable weak solution, satisfying the localized energy inequality, instead, we could construct an approximate system, for which we can deduce Caccioppoli-type inequalities to obtain ``approximate singular set", and then by passing to a limit in an appropriate sense,  we can show that there exists a possible singular set for the limiting weak solution, whose Hausdorff dimension  is at most two. When we  try to apply the similar idea to the full 3D non-stationary system defined on $\Bbb R^3$, we have difficulty in constructing a sequence of the approximate weak solutions, the compactness of which is strong enough to pass to the limit.   Therefore, we leave the proof of partial regularity of the full 3D non-stationary system as an open problem. 
We now formulate our problem more precisely, and state our main result. 
\\

We concentrate on  the following $3D$ valued Hall-MHD system in 
$Q= \R^2 \times (0,T)$.
\begin{align}
\partial _t \bu +  (\bu\cdot \nabla )\bu - \Delta \bu &= -\nabla p +
 (\nabla \times \bB) \times  \bB + \bbf,
\label{1.2}
\\
\partial _t \bB + \nabla \times (\bB\times \bu) - \Delta \bB &= 
-\nabla \times ((\nabla \times \bB) \times  \bB)+  \nabla \times \bg,
\label{1.3}
\\
\nabla \cdot \bu &=0,\quad  \nabla \cdot \bB =0
\label{1.1}  
  \end{align}
together with the initial condition
\be\label{1.3a}
\bu= \bu_0,\quad \bB= \bB_0 \quad \mbox{on}\quad \R^2\times \{0\},
\ee
which satisfy 
\be\label{1.3b}
\nabla \cdot \bu_0 =\nabla \cdot \bB_0 =0 \quad \mbox{on}\quad \R^2.
\ee
Here, $\bu=(u_1, u_2, u_3), \bB=(B_1, B_2, B_3)$, where $u_j=u_j (x_1, x_2, t), B_j=B_j(x_1,x_2, t), j=1,2,3$, and $p=p(x_1,x_2,t)$, $(x,t)=(x_1,x_2,t)\in Q$.
Note that we set $\nu=\mu=1$ for convenience. For the  definition of weak solution see Definition\,\ref{def2.1} below.
The aim of the present paper is to prove the existence of a weak solution to the Hall-MHD 
system \eqref{1.2}--\eqref{1.1},  which is H\"older continuous outside of a possible singular set together with the estimation of its Hausdorff dimension. 
We set
$ L^2_{\rm div} =\{ \bu  \in  L^2 \, |\, \nabla \cdot \bu=0\}$,
where the derivative is defined in the sense of distribution.
We also define $V^2(Q) = L^\infty (0,T; L^2) \cap  L^2(0,T; W^{ 1,\, 2})$. By $V^2_{\rm div}(Q) $ 
we denote the space of all $\bu \in V^2(Q)$ such that $\nabla \cdot \bu =0$ in the sense of distribution in $Q$.  

\vspace{0.2cm}  
\hspace{0.5cm}
Notice that using the formula $(\bu\cdot \nabla)\bu  = (\nabla \times \bu )\times \bu+ \frac {1} {2} |\bu|^2$, one can 
rewrite   \eqref{1.2} into
\be\label{1.4}
 \partial _t \bu + (\nabla \times \bu )\times \bu - \Delta \bu = 
-\nabla \Big(p + \frac {|\bu|^2} {2}\Big)+ (\nabla \times \bB) \times  \bB + \bbf \quad \mbox{in}\quad Q.
\ee
Applying $\nabla \times $ to the both sides of the above, we get
\be\label{1.5}
\partial _t \bfomega + \nabla \times (\bfomega\times \bu) - \Delta \bfomega = 
\nabla \times ((\nabla \times \bB) \times  \bB) + \nabla \times \bbf\quad \mbox{in}\quad Q,
\ee
where $\bfomega $  stands for the vorticity $\nabla \times \bu$. 
Taking the sum of \eqref{1.3} and \eqref{1.5},  we are led to 
\be\label{1.6}
\partial _t \bV + \nabla \times (\bV\times \bu) - \Delta \bV =  \nabla \times (\bbf + \bg)  \quad \mbox{in}\quad Q,
\ee
where
\be\label{1.7}
\bV = \bB+\bfomega. 
\ee
Since $\nabla \cdot \bV=0$, there exists a solenoidal potential $\bv$ such that $\nabla \times \bv= \bV$. 
From \eqref{1.6} we deduce that $\bv$ solves the following system in $Q$,
\begin{align}
 \nabla \cdot  \bv &=0,
\label{1.8}
\\
 \partial _t \bv + (\bv \cdot \nabla) \bv - \Delta \bv &= 
- \nabla \pi + (\nabla \times \bv) \times \bb + \bbf+ \bg,
\label{1.9}
\end{align}
where $\bb= \bv -\bu$.  Clearly, $\nabla \times \bb = \bB$. 

\hspace{0.5cm}
We now introduce the notion of a weak solution to the system \eqref{1.2}--\eqref{1.3b}.  

 \begin{defin}
 \label{def2.1}
Let $\bbf, \bg\in  L^2(Q)$.  We say $(\bu, p, \bB) \in  
V^2_{\rm div}(Q)\times L^2(0,T; L^{2}_{\rm loc})\times V^2_{\rm div}(Q) $  is a 
{\it weak solution} to \eqref{1.1}--\eqref{1.3a} if  
\begin{align}
  &\intl_{Q} (- \bu \cdot \partial _t \bfvarphi + 
\nabla \bu : \nabla \bfvarphi - \bu \otimes \bu : \nabla \bfvarphi ) dxdt 
\cr
&\quad = \intl_{Q} p \nabla \cdot  \bfvarphi  dxdt  + 
\intl_{Q} ((\nabla \times \bB) \times \bB+ \bbf) \cdot \bfvarphi  dxdt  + \intl_{\R^2} \bu_0\cdot \bfvarphi(0)  dx  ,   
\label{2.1}
\\
  &\intl_{Q}( \bB \cdot \partial _t \bfvarphi + \nabla \bB : \nabla \bfvarphi + 
\bB \times \bu : \nabla \times \bfvarphi ) dx dt
\cr
&\quad = \intl_{Q}( (\nabla \times \bB) \times \bB + \bg)\cdot \nabla \times \bfvarphi  dx dt 
+ \intl_{\R^2} \bB_0 \cdot \bfvarphi(0)  dx
\label{2.2}
\end{align} 
for all $\bfvarphi \in  C^\infty _{\rm c}(\R^2\times [0, T)) $. Here 
we used the notation $\bA: \bB  =\sum_{i,j=1}^{3} A_{ ij} B_{ ij}$ 
for matrices $ \bA, \bB \in \R^{3\times 3}$.

\end{defin}

\hspace{0.5cm}
By $ {\mathcal M}^{2,\lambda }_{ \rm loc} (Q) $ we denote the local Morrey space, which is defined in  Section\,3 below.  Our main result is the following theorem. 

\begin{thm}
\label{thm1.1} 
Let $\bu_0 \in L^2_{\rm div}, \bB_0\in L^2$ and  
$\bbf, \bg \in L^2(Q)$. Moreover,  we suppose that $ \bg
\in   {\mathcal M}^{2,\lambda }_{ \rm loc} (Q) $  for some $2<\lambda < 4$. Then,  there exists 
a weak solution $(\bu, p, \bB)\in  V^2_{\rm div}(Q) \times L^2(0,T; L^2_{\rm loc} ) \times V^2_{\rm div}(Q)$ 
of \eqref{1.2}--\eqref{1.3b} being  $\alpha$-H\"older continuous outside of a closed subset set $\Sigma (\bB)\subset Q$  
of Hausdorff dimension less than or equal to two, where $0<\alpha < \frac {\lambda -2} { 2}$. 
\end{thm}

\hspace{0.5cm}
The paper is organized as follows. 
In Section\,2 we discuss local estimates of  weak solutions to the approximate system related to \eqref{1.3}  involving the magnetic field $ \bB$. 
Thanks to the validity of the local energy equality (see \eqref{3.0-4} below) we are able to establish a Caccioppoli-type inequality, which plays a central role  in the  proof of  the  fundamental estimate in Section\,3  (cf. Lemma\,\ref{4.2}).  To achieve this result we make use of an indirect argument together with the fundamental estimate which holds true for the corresponding linear limit system (cf. Lemma\,\ref{4.1}).  
The aim of section Section\,4 is the construction  of an approximate solution to system \eqref{1.2}--\eqref{1.3b} along with the required {\it a priori} estimates. 
Furthermore,  passing to the limit in the approximate system we get  a weak solution to \eqref{1.2}--\eqref{1.3b}. In Section\,5 we prove  that the weak solution constructed in Section\,4 fulfills the required partial regularity property stated  in Theorem\,\ref{thm1.1}, the main result of the paper.  We wish to remark that even for the weak solution to the system under  consideration constructed in a suitable way, a corresponding local energy inequality similar to the case of the Navier-Stokes equation may not be available. For this reason in the proof of the main theorem we are only able to work on the  approximate solutions using Lemma\,\ref{4.2}.  The estimation of the parabolic Hausdorff dimension of the singular set is obtained by Theorem\,\ref{thm6.1}, the proof of which can be found at the end of Section\,5. For readers convenience we  added an appendix which contains the definition of the parabolic H\"older space $ C^{ \alpha , \alpha /2}(Q)$, the parabolic version of the Poincar\'e inequality and an algebraic lemma which will be used  in the proof of Theorem\,\ref{thm6.1}.

\section{Caccioppoli-type inequality for the approximate $\bB$ system}
\label{sec:-2}
\setcounter{secnum}{\value{section} \setcounter{equation}{0}
\renewcommand{\theequation}{\mbox{\arabic{secnum}.\arabic{equation}}}}

Let  $\bg,  \bu \in L^2 (Q) $  be given.   For fixed   $0<\delta <1$ we consider the following system for $\bB$  approximating   \eqref{1.3}
\begin{align}
 &\partial _t \bB - \Delta \bB  
\cr
&= - \nabla \times \Big (\nabla \times \bB  \times \frac {\bB} {1+ \delta |\bB|}\Big) + 
\nabla \times \Big (\bu\times \frac {\bB} {1+ \delta |\bB|}\Big) 
+\nabla \times \bg \quad \mbox{in} \quad Q.
\label{3.0}
\end{align}

We start our discussion with the following notions of a weak solution to \eqref{3.0}.
\begin{defin}
\label{def3.1}
A vector field $\bB \in  V^2 (Q)$ is said to be a 
{\it weak solution to }\eqref{3.0} if 
\begin{align}
&\intl_{Q} (- \bB \cdot \partial _t \bfvarphi + \nabla \bB: \nabla \bfvarphi )  dx dt 
\cr
& \quad = 
- \intl_{Q}  (\nabla \times \bB - \bu)  \times \frac {\bB} {1+ \delta |\bB|} \cdot \nabla \times 
\bfvarphi  dx dt + \intl_{Q} \bg \cdot \nabla \times \bfvarphi  dxdt    
\label{3.0-1}
\end{align}
for all $\bfvarphi \in  C^\infty _{\rm c}(\Omega ) $. 
\end{defin}

\begin{rem}
\label{rem3.2}
Let $\bB$ be a weak solution to \eqref{3.0}.  Then, \eqref{3.0-1} 
yields the existence of the distributional time derivative   
$\bB'\in  L^2(0,T; W^{ -1,\, 2})$, determined by the identity 
\begin{align}
&\intl_{\R^2} \langle \bB'(s), \bfpsi  \rangle dx  
+ \intl_{\R^2} \nabla \bB(s):\nabla \bfpsi  dx   
\cr
& =
-\intl_{\R^2} (\nabla \times \bB(s) - \bu (s))\times  \frac {\bB(s)} {1+ \delta |\bB(s)|} \cdot \nabla \times \bfpsi dx  
 + \intl_{\R^2} \bg(s) \cdot \nabla \times \bfpsi dx  
\label{3.0a}
\end{align}
 for all $\bfpsi \in  W^{ 1,\, 2}(\R^2)$ and   for a.\,e. $s\in  (0,T)$.  Inserting 
$\bfpsi(x,s)= \phi(x,s) (\bB(x,s)- \bLambda ) $   into \eqref{3.0a} with $\phi \in  C^\infty _{\rm c}(Q)$  and a constant vector 
$\bLambda \in \R^3$, integrating the result over $(0,t)$ \,$(t\in  (0,T))$ and 
using integrating by parts, we obtain the following 
local energy equality   
\begin{align}
 &  \frac{1}{2} \intl_{\R^2} \phi(t)  |\bB(t)- \bLambda |^2  dx  + 
\intl_{0}^{t} \intl_{Q} \phi  |\nabla \bB|^2 dx ds 
\cr
&\quad =
 \frac {1} {2} \intl_{0}^{t}  \intl_{ \R^{2}} (\partial _t\phi  +\Delta \phi) |\bB- \bLambda |^2   dx  ds   
\cr
& \qquad \quad +\intl_{0}^{t} \intl_{\R^{2}} (\nabla \times \bB- \bu) \times  
\frac {\bB} {1+ \delta |\bB|} \cdot   ((\bB- \bLambda ) \times  \nabla \phi )   dx ds 
\cr
&\qquad \quad  +\intl_{0}^{t} \intl_{\R^{2}} \phi  \bu \times  
\frac {\bB} {1+ \delta |\bB|} \cdot \nabla \times \bB dx ds
\cr
&\qquad\quad+ \intl_{0}^{t}  \intl_{\R^{2}} 
\Big(\phi \bg \cdot  \nabla \times 
\bB - \bg \cdot (\bB- \bLambda)  \times \nabla \phi \Big) dx  ds.
\label{3.0-4}
\end{align}
\end{rem}

\hspace*{1cm}
First, let us fix some notations which is used throughout the present and subsequent sections. 
Let $ X_0 = (x_0, t_0)\in \R^{3}$ and $ 0<r<+\infty$ by $ Q_r= Q_r(X_0)$ we denote the parabolic cylinder $ B_r(x_0) \times (t_0-r^2, t_0)$.  
Furthermore, for a function  $f\in L^1(Q_r) $ we define  
\[
f_{ r, X_0}:= f_{ Q_r}= \intmw_{Q_r} f dxdt = \frac{1}{\mes Q_r} \intl_{Q_r} f  dx dt,\
\]
where $ \mes Q_r$ stands for the three dimensional Lebesgue measure of $ Q_r$. 

\hspace{0.5cm}
Let $ 0<\rho < r$. We call $\theta \in  C^\infty (\R^3)$  a {\it cut-off 
function} suitable  for $Q_{r}$ and $Q_\rho $ if $0\le \theta \le 1$  in $\R^3$, $\theta \equiv 1$ on $Q_{\rho } $, 
$\theta \equiv 0$ in $(\R^3\setminus B_r ) \times (t_0-r^2, t_0)\cup  \R^2\times (-\infty , t_0 - r^2)$  
and $|\partial _t \theta| +|\nabla \theta |^2 + |\nabla ^2 \theta |\le c(r-\rho )^{-2} $ in $\R^3$.

\hspace{0.5cm}
Now, we state the following Caccioppoli-type inequality.

\begin{lem}
\label{lem3.3}
Let $\bg\in  L^2 (Q), \bu \in L^4(Q)$ be given, and let $\bB \in  V^{2} (Q)$ be a weak solution to \eqref{3.0}.  
Then, for every cylinder $\overline{Q}_r= \overline{Q_r(X_0)}  \subset  Q$ and 
$0< \rho < r$ there holds 
\begin{align}
& \esssup_{t\in (t_0- r^2, t_0)} \intl_{B_r} \theta ^{4} |\bB- \bB_{r, X_0} |^2  dx  
+\intl_{Q_r}  \theta ^{4} |\nabla  \bB|^2  dx dt 
\cr
&  \quad \le  \frac {c} {(r-\rho )^2} (1+ |\bB_{r, X_0} |^2)\intl_{Q_r} |\bB- \bB_{r, X_0} |^2  dx dt  
\cr
& \qquad +
\frac {c} {r-\rho }\bigg(\intl_{Q_r}  \theta ^{3+ \gamma  }  |\bB- \bB_{r, X_0} |^4 dx dt\bigg)^{1/2} \bigg( \intl_{Q_r}  \theta ^{3-\gamma } |\nabla  \bB|^2  dx dt\bigg)^{1/2}
\cr
&\qquad + c\intl_{Q_r} ( |\bg|^2 + \theta ^4|\bB|^2 |\bu|^2  )  dx dt
\label{3.5}
\end{align}
 for all  {\it cut-off function } $ \theta $ suitable  for $Q_{r}$ and $Q_\rho $   
$( \gamma \in [-3,3])$, and 
\begin{align}
 E(\rho )^4 &\le 
\frac {cr^4} {(r-\rho )^4} (1+ |\bB_{r, X_0} |^2) (G(r)^4 + F(r)^4) 
\cr
& \qquad +
\frac {cr^6} {(r-\rho)^6 } (G(r)^6 + F(r)^{6}) 
\cr
&\qquad + \frac {c} {(r-\rho )^2} \Bigg\{\intl_{Q_r} |\bg|^2 dxdt + |\bB_{r, X_0}|^2  
\intl_{Q_r}  |\bu|^2  dx dt \Bigg\} (G(r)^2 + F(r)^2)
\cr
&\qquad  +  \frac {cr^4} {(r-\rho )^4} \intl_{Q_r}  |\bu|^4    dx dt  (G(r)^4 + F(r)^4), 
\label{3.6}
\end{align}
where $c=\const>0$ denotes a universal constant, and 
\begin{align*}
E(r)= E(r, X_0) &=   \bigg(\intmw_{Q_r(X_0)} |\bB-  \bB_{r, X_0}  |^4  dxdt\bigg)^{1/4},
\\
F(r)= F(r, X_0) &=   \bigg(r^{-2} \intl_{Q_r(X_0)} |\nabla \bB  |^2  dxdt\bigg)^{1/2},\quad 
\\
G(r)= G(r, X_0) &=   \bigg(\intmw_{Q_r(X_0)} |\bB-  \bB_{r, X_0}  |^2  dxdt\bigg)^{1/2},
0<r< \sqrt{t_0}.
\end{align*}

\end{lem}
{\it Proof }  Let $\overline{Q} _r = \overline{Q_r(X_0)}  \subset  Q$  be a fixed cylinder.  
For $0<\rho < r$  we take a cut-off function 
$\theta \in  C^\infty (\R^3)$  suitable  for $Q_{r}$ and $Q_\rho $.  
 
 \hspace{0.5cm}
From \eqref{3.0-4} with $\phi = \theta ^{4} $ and $\bLambda = \bB_{r, X_0} $  we obtain the following Caccioppoli-type inequality
\begin{align}
& \esssup_{t\in (t_0- r^2, t_0)} \intl_{B_r} \theta ^{4} |\bB- \bB_{r, X_0} |^2  dx  
+\intl_{Q_r}  \theta ^{4} |\nabla  \bB|^2  dx dt 
\cr
&\quad\le \frac {c} {(r-\rho )^2} \intl_{Q_r} |\bB- \bB_{r, X_0} |^2  dxdt  +  
c\intl_{Q_r} |\bg|^2 + \theta ^4|\bB|^2 |\bu|^2  dxdt 
\cr
&\qquad\qquad + \frac {c} {r-\rho }\intl_{Q_r} \theta^3   |\nabla \bB|\, 
|\bB|\,|\bB- \bB_{r, X_0}|dxdt
\cr
&  \quad =  \frac {c} {(r-\rho )^2} \intl_{Q_r} |\bB- \bB_{r, X_0} |^2  dx dt  +   
c\intl_{Q_r} |\bg|^2 + \theta ^4|\bB|^2 |\bu|^2    dx dt +J. 
\label{3.1}
\end{align}
Let $\gamma \in [-3,3]$.  Applying H\"older's and Young's inequality, we estimate
\begin{align*}
J& \le \frac {c} {(r-\rho )^2} |\bB_{r, X_0}   |^2 \intl_{Q_r} |\bB- \bB_{r, X_0} |^2 dx dt
\\
&\qquad\qquad+ \frac {c} {r-\rho }\bigg(\intl_{Q_r}  \theta ^{3+\gamma }  |\bB- \bB_{r, X_0} |^4 dx dt\bigg)^{1/2} \bigg( \intl_{Q_r} \theta ^{3-\gamma}   |\nabla  \bB|^2  dx dt\bigg)^{1/2}
\\
&\qquad \qquad+ \frac {1} {2}\intl_{Q_r}  \theta ^4 |\nabla  \bB|^2  dx dt.
\end{align*}
Inserting the estimate of $J$ into \eqref{3.1}, we are led to 
\begin{align}
& \esssup_{t\in (t_0- r^2, t_0)} \intl_{B_r} \theta ^{4} |\bB- \bB_{r, X_0} |^2  dx  
+\intl_{Q_r}  \theta ^{4} |\nabla  \bB|^2  dx dt 
\cr
&  \quad \le  \frac {c} {(r-\rho )^2} (1+ |\bB_{r, X_0} |^2)\intl_{Q_r} |\bB- \bB_{r, X_0} |^2  dx dt  
\cr
& \qquad +
\frac {c} {r-\rho }\bigg(\intl_{Q_r}  \theta ^{3+\gamma }  |\bB- \bB_{r, X_0} |^4 dx dt\bigg)^{1/2} \bigg( \intl_{Q_r} \theta ^{3-\gamma }  |\nabla  \bB|^2  dx dt\bigg)^{1/2}
\cr
&\qquad + c\intl_{Q_r} ( |\bg|^2 + \theta ^4 |\bB|^2 |\bu|^2 )   dx dt.
\label{3.8}
\end{align}
This proves (\ref{3.5}).  On the other hand, by means of Sobolev's embedding theorem we get 
\begin{align}
 & \intmw_{Q_r} \theta ^4 |\bB- \bB_{r, X_0 } |^ 4 dxdt 
\cr
&\le c r^{-4} \|\theta ^2 (\bB - \bB_{r, X_0}) \|^2_{L^\infty(t_0-r^2, t_0; L^2(B_r)) } 
\| \nabla \bB \|^2_{L^2(Q_r)}  
\cr
&\qquad +  c r^{-4} (r-\rho )^{-2} 
\|\theta ^2 (\bB - \bB_{r, X_0}) \|^2_{L^\infty(t_0-r^2, t_0; L^2(B_r)) }
\|\bB - \bB_{r, X_0} \|^2_{L^2(Q_r) }
\cr
& \le  \frac {c} {(r-\rho )^2}\|\theta ^2 (\bB - \bB_{r, X_0}) \|^2_{L^\infty(t_0-r^2, t_0; L^2(B_r)) } 
(F(r)^2+ G(r)^2). 
\label{3.8a}
\end{align}

Combining \eqref{3.8} with $\gamma =1$ and \eqref{3.8a} with help of Young's inequality,  we get 
\begin{align}
&  \intmw_{Q_r}  \theta ^{4} | \bB- \bB_{r, X_0} |^4  dx dt 
\cr
&  \quad \le  \frac {cr^4} {(r-\rho )^4} (1+ |\bB_{r, X_0} |^2)
G(r)^2(F(r)^2 + G(r)^2) 
\cr
& \qquad + \frac {cr^6} {(r-\rho)^6 } (F(r)^6 + G(r)^6) 
\cr
&\qquad + \frac {c} {(r-\rho )^2}
\intl_{Q_r} ( |\bg|^2 + \theta ^4|\bB|^2 |\bu|^2 ) dxdt (F(r)^2 + G(r)^2).
\label{3.1b}
\end{align}

Estimating $|\bB|^2 \le 2|\bB- \bB_{r, X_0} |^2 + 2 |\bB_{r, X_0} |^2$ and applying Young's inequality, we obtain \eqref{3.6}. Thus, the proof of the Lemma is complete. \hfill \Beweisende


\begin{rem}
\label{rem3.4}
From \eqref{3.5} with $\gamma = -1$  along with Young's inequality we get  
\begin{align}
&\bigg(  \frac {1} {r^2}  \esssup_{t\in (t_0-\rho ^2, t_0)} \intl_{B_\rho } 
|\bB(t)- \bB_{r, X_0} |^2 dx   \bigg)^{1/2} + F(\rho )  
\cr
& \quad \le 
\frac {cr} {r-\rho } \Big\{(1+|\bB_{r, X_0} |) E(r) + E(r)^2\Big\} 
\cr
 & \qquad \quad +  \frac {c} {\rho } \Big\{\|\bu\|_{2, Q_r} (E(r) + 
|\bB_{r, X_0}|)+  \|\bg   \|_{2, Q_r}\Big\}. 
\label{3.4-1}
\end{align}

Furthermore, using the parabolic Poincar\'{e}-type inequality 
(cf. Lemma\,\ref{lemA.1}, appendix below), 
we find 
\begin{align}
& \intmw_{Q_r} |\bB- \bB_{r, X_0} |^2 dxdt 
\cr
&\quad \le c(1+ |\bB_{r, X_0} |^2) r^{-2} \intl_{Q_r} |\nabla \bB|^2  dx dt      
\cr
& \qquad  + c (1+  |\bB_{r, X_0} |^2)  r^{-2}\intl_{Q_r} (|\bg|^2 +|\bu|^2) dx dt 
\cr
 & \qquad + C_1 r^{-2}  \intl_{Q_r}( |\nabla \bB|^2 + |\bu|^2 ) dxdt  
\intmw_{Q_r}   |\bB- \bB_{r, X_0} |^2  dxdt   
\label{3.10}
\end{align}
with an absolute constant $C_1>0$.  Thus, assuming that  
\be\label{3.11}
C_1 \Bigg\{r^{-2}  \intl_{Q_r} |\nabla \bB|^2 dxdt +  4 \bigg(\intl_{Q_r} |\bu|^4 dxdt\bigg)^{1/2}\Bigg\}   
\le \frac {1} {2},
\ee
\eqref{3.10}  leads to 
\begin{align}
 G(r) & \le c (1+ |\bB_{r, X_0} |) (F(r)+ H(r)),
\label{3.13}
\end{align}
where
\[
H(r)= H(r, X_0) = r^{-1} \|\bg\|_{2, Q_r} + \|\bu\|_{4, Q_r},\quad 
0<r < \sqrt{t_0}.     
\]
Substituting $G(r)$ on the right of \eqref{3.6}   by \eqref{3.13}, setting  $\rho =\frac {r} {2}$ therein, we arrive   at 
\begin{align}
 E(r/2) \le C_2 (1+ |\bB_{r, X_0}|^2 )  \Big\{F(r) + F(r)^2 + H(r)+ H(r)^2\Big\}
\label{3.14}
\end{align}
with an absolute constant $C_2>0$, provided \eqref{3.11} is fulfilled.  

\hspace*{0.5cm}
From \eqref{3.4-1} with $\rho = \frac {r} {2}$ we deduce that 
\begin{align}
 F(r/2) \le C_3 (1+ |\bB_{r, X_0}| )  \Big\{E(r) + E(r)^2 + H(r)+ H(r)^2\Big\}
\label{3.15}
\end{align}
with an absolute constant $C_3>0$. 

\end{rem}

%
%
\section{Blow-up lemma} 
\setcounter{secnum}{\value {section}
\setcounter{equation}{0}
\renewcommand{\theequation}{\mbox{\arabic{secnum}.\arabic{equation}}}} 

In what follows we define the space 
\[
V^2(Q_r)= L^\infty(t_0-r^2, t_0; W^{1,\, 2}(B_r(x_0)))\cap L^\infty(t_0-r^2, t_0; L^2(B_r(x_0)))
\]
for $ X_0=(x_0, t_0)$ and $ 0<r<+\infty$. 

\vspace{0.3cm}
We begin our discussion with the following fundamental estimate for solutions to the model problem in $ Q_1=Q_1(0,0)$, which will be used in the blow-up lemma below.
\vspace{0.3cm}
\begin{lem}
\label{lem4.1}
Let $\bLambda \in \R^3$. Let $\bW\in L^4(Q_1)$ such that $ \bW|_{ Q_\sigma } 
\in  V^{2} (Q_\sigma )$ for all $ 0<\sigma <1$ solves 
\be\label{4.1}
\partial _t \bW-\Delta \bW = - \nabla \times ((\nabla \times \bW) \times  \bLambda) 
\quad \mbox{in}\quad Q_1 
\ee
in sense of distributions, i.\,e.   
\begin{align}
&\intl_{B_1} \bW(t) \cdot \bPhi  (t) dx + \intl_{-1}^{t}  \intl_{B_1} (- \bW\cdot \partial _t \bPhi + 
\nabla \bW : \nabla \bPhi )   dx ds
\cr
&\quad =   - \intl_{-1}^{t} \intl_{B_1} ((\nabla \times \bW) \times \bLambda) \cdot  
\nabla \times \bPhi    dxds
\label{4.2}
\end{align}
for all $\bPhi \in  W^{ 1,\, 2}(Q_1)$ compactly supported in  $Q_1$, for a.\,e. 
$t\in (-1,0)$.   
Then, 
\be\label{4.3}
\bigg(\intmw_{Q_\tau } |\bW- \bW_{Q_\tau }|^4 dxdt  \bigg)^{1/4} \le 
C_0 \tau (1+ |\bLambda  |^{5} ) 
\bigg(\intmw_{Q_1} |\bW- \bW_{Q_1 }|^4 dxdt\bigg)^{1/4}
\ee
for all $0<\tau <1$, where $C_0>0$ denotes a universal constant. 
\end{lem}

{\it Proof}  Since the assertion is trivial for $\frac {1} {4} < \tau < 1$, we may assume that $0< \tau \le  \frac {1} {4}$. 
Let $\zeta \in  C^\infty _{\rm c} (\R^3)$  be a suitable cut-off function for $Q_\tau $ and 
$Q_{1/2} $.   Inserting  the admissible test function  
$\bPhi = \zeta ^{2m} (\bW-\bW_{B_1}) $\, 
$(m\in \N)$ into \eqref{4.2}, by using Cauchy-Schwarz's inequality along with 
Young's inequality, we are led to  
\begin{align}
&\esssup_{t\in (-1, 0)} \intl_{B_1} \zeta ^{2m} |\bW(t)|^2  dx +   
\intl_{Q_1} \zeta ^{2m} |\nabla \bW|^2 dxdt 
\cr
&\qquad \qquad \le c (1+ |\bLambda |^2) \intl_{Q_1} \zeta ^{2m-2} |\bW- \bW_{Q_1} |^2 dxdt.     
\label{4.4}
\end{align}
If $\bW$ is smooth in $Q_1$, since \eqref{4.1} is  a linear system,  the same inequality holds true for $D^\alpha \bW $ in place of $\bW$ for  any multi-index $\alpha $. 
By a standard mollifying  argument together with Sobolev's embedding theorem 
we see that $\bW$ is smooth in $Q_1$. By an iterative application of \eqref{4.4} 
with $ m=4,3,2,1$ we obtain   
\be\label{4.5}
\esssup_{t\in (-1,0)} \intl_{B_1} \zeta ^{8} |D^\alpha \bW|^2 dx \le 
c (1+ |\bLambda |^8) \intl_{Q_1} |\bW- \bW_{Q_1} |^2 dxdt\quad\forall\, |\alpha |\le 3.     
\ee  
By means of Sobolev's embedding theorem and Jensen's inequality we get 
\be\label{4.6}
 \|\nabla \bW\|_{\infty , Q_{1/2} }^4 \le   c (1+ |\bLambda |^{16} ) \intl_{Q_1} |\bW- \bW_{Q_1} |^4 dxdt.
\ee
Applying Poincar\'e's inequality, we arrive at
\be\label{4.7}
\intmw_{Q_\tau } |\bW- \bW_{Q_\tau } |^4 dx dt\le c\tau^4  (1+ |\Lambda |^4) \|\nabla \bW\|_{\infty , Q_{1/2} }^4.  
\ee
Combination  of  \eqref{4.6} and \eqref{4.7} gives the desired estimate.
\hfill \Beweisende 
\vspace{0.3cm}
\hspace*{1cm}

\hspace{0.5cm}
In our discussion below we make use of the notion of the Morrey space. 
Let $ K \subset Q$ be a compact set. Define, $ d_K = \min \{t\in (0,T) \,|\,t\in K \}$. We say $f$ belongs to the Morrey space ${\cal M}^{p, \lambda } (K)$ if
\[
[f]_{{\cal M}^{p, \lambda }, K} := \sup\bigg\{r^{-\lambda } \intl_{Q_r(X_0)} 
|f|^{p}  dxdt  \,\bigg|\, X_0\in K, 0<r\le 
d_K\bigg\}   <+\infty.  
\] 
Furthermore, by  $ f\in  {\cal M}^{p, \lambda }_{ \rm loc} (Q)$ we mean 
 $ f\in {\cal M}^{p, \lambda } (K)$ for all compact set $ K \subset Q $. 

 \vspace{0.2cm}   
 \hspace{0.5cm}
 Now we are ready to state the following key lemma. 

\begin{lem}
\label{lem4.2}
Let $\bg \in {\cal M}^{2, \lambda}_{\rm loc} (Q)$  for some 
$2< \lambda < 4$. 
For every  $0<\tau < \frac {1} {2},  0<M, L < +\infty$, compact set  $K \subset Q $ and
$0<\alpha < \frac {\lambda -2} {2}$,  
there exist positive numbers 
 $\var _0= \var _0(\tau , M, L, K, \alpha )$, $R _0 = R _0(\tau , M, L, K, \alpha ) < d_K$ and $\delta _0 = \delta _0 (\tau , M, L, K, \alpha ) \le 1$ 
such that, if  $\bB \in V^{2} (Q)$  is a 
weak solution to \eqref{3.0}with $0< \delta \le \delta _0$ and  $\bu \in L^{8/(4-\lambda ) } (Q)$  such that 
\be\label{4.8a}
\|\bu \|_{8/(4-\lambda ), Q } \le L,
\ee
and if for $X_0\in K$  and  $0<R\le R_0$  the following condition is fulfilled
\be\label{4.8}
|\bB_{R, X_0 }| \le M, \quad E (R, X_0) + R^{\alpha } \le \var _0,
\ee 
then there holds
\be\label{4.9}
E (\tau R, X_0) \le 2 \tau C_0 (1+M^{5}) (E (R, X_0) + R^{\alpha}), 
\ee
where $C_0>0$   stands for the constant appearing on the right hand side of \eqref{4.3}. 
\end{lem}

{\it Proof }  Assume the assertion of the Lemma is not true.  Then there exist 
$0< \tau < \frac {1} {2}, 
0<M, L < +\infty$, a compact set $K \subset  \Omega $  and $0<\alpha < \frac {\lambda -2} {2}$  
as well as sequences $ \{\varepsilon _k\}, \{\delta _k\}\subset  (0,1)$ with $ \varepsilon _k, \delta _k \rightarrow 0$ as $ k \rightarrow +\infty$, 
 $\{R_k\} \subset  (0, d_K)$, $\{X _k\} =\{(x_k,t_k)\}\subset  K$,
$\{ \bu^{ (k)}\} \subset  L^{ 8/(4-\lambda )}(Q)$ fulfilling 
\begin{equation}
\| \bu ^{ (k)}\|_{ 8/(4-\lambda )} \le L\quad  \forall\,k\in \N,
\label{4.8b}
\end{equation}
and a sequence $\{\bB^{(k)}\} \subset  V^{2} (Q)$,  being a weak solutions  to \eqref{3.0} 
replacing $ \bu $ by $ \bu ^{ (k)}$   and  $ \delta $ by  $ \delta _k$ respectively, 
such that 
\be\label{4.10}
|\bB^{(k)} _{ R_k, X_k}| \le M, \quad E_k (R_k, X_k) + R_k^{\alpha } = \var _k
\ee 
and 
\be\label{4.11}
E_k (\tau R_k, X_k) > 2 \tau C_0 (1+M^{5}) (E_k (R_k, X_k)+ R_k^{\alpha } ).
\ee
Here we have used the notation
\[
E_k (r, X_k) = \bigg(\intmw_{Q_r(X_k)} |\bB^{(k)} - \bB^{(k)}_{r, X_k}  |^4 dxdt\bigg)^{1/4},\quad 
X_k \in K, 0<r\le d_K  
\]
$(k\in \N)$. Note that \eqref{4.10} yields $R_k \rightarrow 0$ as $k \rightarrow +\infty $. 

\hspace*{0.5cm}
Next, for $Y:=(y,s)\in Q_1(0) $ we define
\begin{align*}
\bW_k (Y) &= \frac {1} {\var _k} (\bB^{(k)} (x_k + R_k y, t_k+R_k^2 s) - \bB^{(k)} _{ R_k, X_k}),
\\     
\bv_k (Y) &=  \bu^{(k)} (x_k + R_k y, t_k+R_k^2 s), 
\\     
\bg_k (Y) &=  \bg(x_k + R_k y, t_k +R_k^2 s),
\end{align*}
$( k\in \N)$.  Furthermore, we set
\[
\mathscr{E}_k (\sigma ) = \bigg(\intmw_{Q_\sigma } |\bW_k- (\bW_k)_{Q_\sigma }|^4  dyds\bigg)^{1/4},
\quad 0<\sigma \le 1.   
\]
Then \eqref{4.10} and \eqref{4.11} turn into 
\be\label{4.10-1}
|\bB^{(k)} _{ R_k, X_k}| \le M, \quad \mathscr{E}_k (1) + \frac { R_k ^{\alpha  }} {\var _k}= 1,
\ee 
and 
\be\label{4.11-1}
\mathscr{E}_k (\tau ) > 2 \tau C_0 (1+M^{5}) \Big(\mathscr{E}_k (1)+ 
\frac { R_k^{\alpha  }} {\var _k}\Big)= 2 \tau C_0 (1+M^{5})
\ee
respectively. 

\hspace*{1cm}
Using the chain rule,  restriction of system \eqref{3.0} to $ Q_{ R_k}(X_k)$ 
takes the form 
\begin{align}
& \partial _t \bW_k- \Delta \bW_k 
\cr
& \quad = - 
\nabla \times \bigg((\nabla \times  \bW_k) \times  
\frac {\var _k \bW_k + \bB^{(k)}_{R_k, X_k}} {1+ \delta _k |\var _k \bW_k + \bB^{(k)}_{R_k, X_k}|}\bigg)
\cr
& \qquad 
+ \frac {R_k} {\var _k} \nabla \times   \bigg(\bv_k \times  \frac {\var _k \bW_k + \bB^{(k)}_{R_k, X_k}} {1+ \delta _k |\var _k \bW_k + \bB^{(k)}_{R_k, X_k}|}\bigg) 
+\frac {R_k} {\var _k} \nabla \times  \bg_k
\label{4.12}
\end{align}
in $Q_1$. Thus, $\bW_k \in V^{ 2}(Q_1)$ is a weak solution to \eqref{4.12}.  

\hspace*{0.5cm}
Let $0<\sigma <1$. Using the transformation formula, noticing that $|\bB^{(k)}_{R_k, X_k}|\le M  $,  
the Caccioppoli-type inequality 
\eqref{3.4-1} with $r= R_k$ and $\rho = \sigma  R_k$  turns into 
\begin{align}
 & \|\bW_k\|_{L^\infty (-\sigma ^2,0; L^2(B_\sigma ))  } + 
 \|\nabla \bW_k\|_{2, B_\sigma } 
\cr
&\quad \le c (1-\sigma )^{-1} 
\Big((1+M) \mathscr{E}_k(1) + 
\var_k  \mathscr{E}_k(1)^2\Big)  
\cr
&\qquad + \frac {c R_k^{-1}} {\var _k  } 
\Big(\|\bu^{(k)}  \|_{2, Q_{R_k}(X_k)} (\var _k \mathscr{E}_k(1)+ M) + \|\bg \|_{2, Q_{R_k}(X_k)} \Big).
\label{4.12-1}
\end{align}
As $\bg_k \in  {\cal M}^{2, \lambda }( K) $ observing \eqref{4.10-1}, we see that 
\be\label{4.12-1a}
\frac {R_k^{-1} } {\var _k}\|\bg \|_{2, Q_{R_k}(x_k)}  \le \frac { R_k^{(\lambda-2)/2}} {\var _k} 
[\bg]_{ {\cal M}^{2, \lambda }( K) } 
\le  R_k^{(\lambda -2)/2-\alpha }   [\bg]_{ {\cal M}^{2, \lambda }(K) }.
\ee
Similarly, by \eqref{4.8b} and \eqref{4.10-1} we get 
\be\label{4.12-1b}
\frac {R_k^{-1} } {\var _k}\|\bu^{(k)}  \|_{2, Q_{R_k}(x_k)}  \le c\frac { R_k^{(\lambda-2)/2}} {\var _k} 
\|\bu^{(k)} \|_{8/4-\lambda, Q}   
\le  cR_k^{(\lambda -2)/2-\alpha }  L. 
\ee
Thus, from \eqref{4.12-1} with help of \eqref{4.12-1a}, \eqref{4.12-1b}  and \eqref{4.10-1} we obtain 
\be\label{4.12-2}
\|\bW_k\|_{L^\infty (-\sigma ^2,0; L^2(B_\sigma ))  }+ \|\nabla \bW_k\|_{2, Q_\sigma } \le c(1-\sigma )^{-1} (M+ 1)  + 
c ([\bg]_{ {\cal M}^{2, \lambda }( K )}+ L).   
\ee
In addition, in view of  \eqref{4.10-1} we estimate 
\be\label{4.12-3}
\|\bW_k\|_{4, Q_1} = (\mes B_1)^{1/4}  \mathscr{E}_k (1) \le (\mes B_1)^{1/4}.  
\ee
From \eqref{4.12-2} and \eqref{4.12-3} it follows that $\{\bW_k\}$ is bounded in $V^{2}(Q_\sigma )$ for all $0<\sigma <1$ and bounded in $L^4(Q_1)$.  
Thus, by means of reflexivity,  eventually passing to subsequences, 
we get $\bW\in  L^4(Q_1)$ with $ \bW \in V^{2}(Q_\sigma )$ for all $ 0<\sigma <1$ and  $\bLambda \in \R^3$ such that 
\begin{align}
 \bW_k &\rightarrow \bW \quad \mbox{{\it weakly in}} \quad L^4(Q_1) \quad \mbox{as} \quad k \rightarrow +\infty, 
\label{4.13}
\\
\nabla \bW_k &\rightarrow \nabla \bW \quad \mbox{{\it weakly in}} \quad L^{2}(Q_\sigma )
\quad \mbox{as} \quad k \rightarrow +\infty\quad \forall\, 0<\sigma <1, 
\label{4.13-1}
\\
\bW_k &\rightarrow \bW \quad \mbox{{\it weakly$^\ast$ in}} \quad L^\infty (-\sigma ^2, 0; L^2(B_\sigma )) \quad \mbox{as} \quad k \rightarrow +\infty\quad \forall\, 0<\sigma <1, 
\label{4.13-1b}
\\
\bLambda_k &\rightarrow \bLambda  \quad\mbox{in}\quad \R^3 \quad 
\mbox{as}\quad k \rightarrow +\infty. 
\label{4.14}
\end{align}
On the other hand, from \eqref{4.12} we deduce that the sequence of distributive time derivative $\{\bW'_k\} $  
is bounded in $L^{4/3}(-\sigma ^2, 0; W^{ -1,\, 4/3}(B_\sigma))$. From this fact  together with \eqref{4.13}  it follows that
\be\label{4.19}
\bW_k \rightarrow \bW \quad \mbox{{\it strongly in}} \quad L^{2}(Q_\sigma ) \quad \mbox{as} \quad k \rightarrow +\infty \quad \forall\, 0<\sigma < 1.
\ee   
Thus, we are in a position to carry out the passage to the limit $k \rightarrow+\infty $ in the weak formulation of \eqref{4.12} to deduce that $\bW$ is a weak solution to the linear system \eqref{4.1}. 
\hspace*{0.5cm}
Our next aim is to prove the strong convergence of $\bW_k \rightarrow \bW$  
in $L^4(Q_\sigma )$ \, $(0<\sigma <1)$.
We first state the following energy equality, 
\begin{align}
& \frac {1} {2} \intl_{B_1} \phi^2(t) |\bW_k(t)|^2  dy + 
 \intl_{-1}^{t}  \intl_{B_1}\phi ^2|\nabla \bW_k|^2  dyds     
\cr
& = \frac {1} {2}  \intl_{-1}^{t}  \intl_{B_1} (\partial _t \phi^2 + \Delta \phi^2 ) |\bW_k|^2 dyds
\cr
& \qquad + \intl_{-1}^{t}  \intl_{B_1} (\nabla \times \bW_k) \times  \frac {\var _k \bW_k + \bB^{(k)}_{R_k, X_k}} {1+ \delta _k |\var _k \bW_k + \bB^{(k)}_{R_k, X_k}|}  \cdot 
(\bW_k \times \nabla \phi^2) dyds 
\cr
& \qquad + \frac {R_k} {\var _k} 
 \intl_{-1}^{t}  \intl_{B_1}\bigg\{\bv_k \times  \frac {\var _k \bW_k + \bB^{(k)}_{R_k, X_k}} {1+ \delta _k |\var _k \bW_k + \bB^{(k)}_{R_k, X_k}|} + \bg_k \bigg\}
\nabla \times (\phi ^2\bW_k) dyds
\label{4.19a}
\end{align}
for all $t\in  [-1,0]$. 
 In view of \eqref{4.13}, \eqref{4.13-1}, \eqref{4.14} and \eqref{4.19}  on both sides of \eqref{4.19a} 
with $t=0$ letting $k \rightarrow +\infty $, we infer  
\begin{align}
&\lim_{k \to \infty } \bigg(\frac {1} {2} \intl_{B_1} \phi^2(0) |\bW_k(0)|^2  dy + 
\intl_{Q_1} \phi ^2|\nabla \bW_k|^2  dyds  \bigg)     
\cr    
& = \frac {1} {2} \intl_{Q_1} (\partial _t \phi^2 + \Delta \phi^2 ) |\bW|^2 dyds
- \intl_{Q_1} (\nabla \times \bW)  \times \bLambda \cdot (\bW \times \nabla \phi^2) dyds. 
\label{4.19b}
\end{align}
Since $\bW$ is a weak solution to \eqref{4.1},  there holds 
\begin{align}
&\frac {1} {2} \intl_{B_1} \phi^2(0) |\bW(0)|^2  dy + 
\intl_{Q_1} \phi ^2 |\nabla \bW|^2  dyds       
\cr  
& = \frac {1} {2} \intl_{Q_1} (\partial _t \phi^2 + \Delta \phi^2 ) |\bW|^2 dyds
-\intl_{Q_1} (\nabla \times \bW)  \times \bLambda \cdot (\bW \times \nabla \phi^2) dyds. 
\label{4.19c}
\end{align}
Noticing that 
\[
\begin{cases}
(\phi(0) \bW_k(0), \phi \nabla \bW_k) \rightarrow 
(\phi(0) \bW(0), \phi \nabla \bW) \quad 
\\[0.2cm]
\mbox{{\it weakly in} }
\quad L^2(B_1) \times L^2(Q_1) \quad 
\mbox{as}\quad k \rightarrow +\infty 
\end{cases}
\]
from \eqref{4.19b} and \eqref{4.19c},  we deduce that  
\[
\nabla \bW_k \rightarrow \nabla \bW \quad \mbox{{\it strongly in}} \quad 
L^2(Q_\sigma ) \quad \mbox{as} \quad k \rightarrow +\infty
\quad  \forall\,0<\sigma <1.
\]
Accordingly,
\be\label{4.20}
\lim_{k \to \infty } \mathscr{E}_k (\sigma ) = \mathscr{E} (\sigma )\quad \forall\,  0<\sigma <1, \ee where $\mathscr{E}(\sigma )= {\D \bigg(\intmw_{B_\sigma } |\bW- \bW_{B_\sigma }  |^4 dy\bigg)^{1/2}}$. 
In particular, thanks to \eqref{4.20} (with $\sigma =\tau $) from \eqref{4.11-1} we get 
\be\label{4.20-1}
\mathscr{E} (\tau ) \ge 2 \tau C_0 (1+ M^5). 
\ee 
 
Since $\bW$ is a weak solution to \eqref{4.1} and $|\bLambda |\le M$,  appealing to  Lemma\,\ref{lem4.1}, we find 
\be\label{4.22}
\mathscr{E}(\tau ) \le \tau C_0(1+ M^5) \mathscr{E}(1).    
\ee
On the other hand,  by virtue of   the lower semi continuity of the norm 
together with   \eqref{4.11-1} and \eqref{4.20} we get 
\begin{align*}
\mathscr{E}(1) &\le  \liminf_{k \to \infty } \Big(\mathscr{E}_k(1) + 
\frac {R_k^{\alpha } } {\var _k}\Big) \le    \frac {1} { 2\tau C_0(1+ M^5)} \lim_{k \to \infty } \mathscr{E}_k(\tau )
\\     
&= \frac {1} { 2\tau C_0(1+ M^5)} \mathscr{E}(\tau ).
\end{align*}
Estimating the right of \eqref{4.22} by the inequality,  we have just obtained we are led to 
$\mathscr{E}(\tau ) \le \frac {1} {2} \mathscr{E}(\tau )$ and hence  
$\mathscr{E}(\tau ) =0$, which contradicts to \eqref{4.20-1}. Whence, the assumption 
cannot be true, which completes the proof of the Lemma. \hfill \Beweisende 

%
%
\section{Construction of approximate solutions} 
\setcounter{secnum}{\value {section}
\setcounter{equation}{0}
\renewcommand{\theequation}{\mbox{\arabic{secnum}.\arabic{equation}}}}

The aim of the present section is to construct a weak solution of the Hall-MHD system 
\eqref{1.2}--\eqref{1.3b} as a limit of a sequence of solutions to the a corresponding approximate system. 
As we will see in the following section, such solution will satisfy the desired partial regularity 
as stated in  the main result of the present  paper. 

\hspace*{0.5cm}
Let $\{\delta _m\} \subset   (0,1)$ \, $(m\in \N)$ be a sequence, such that $\delta _m \rightarrow 0$ as $m \rightarrow +\infty $.  Now, we consider the following approximate system  
\begin{align}
\partial _t \bu_m +  &\frac {\bfomega _m } {1+\delta _m|\bB_m|} \times \bu_m - \Delta \bu_m 
\cr
&\qquad \qquad = -\nabla p_m +
 (\nabla \times \bB_m) \times  \frac {\bB_m } {1+ \delta_m |\bB_m|} + \bbf,
\label{5.2}
\\
\partial _t \bB_m+ &\nabla \times \left(\frac {\bB _m } {1+\delta _m|\bB_m|} \times \bu_m \right)-  \Delta \bB_m 
\cr
&\qquad \qquad = 
-\nabla \times \Big (\nabla \times \bB_m \times  \frac {\bB_m } {1+ \delta_m |\bB_m|} \Big)+  \nabla \times \bg,
\label{5.3}\\
&\qquad \qquad\nabla \cdot \bu_m =0, \quad \nabla \cdot \bB_m=0,
\label{5.1}  
\end{align}
in $Q=\R^2\times  (0,T)$, 
together with the initial condition 
\be\label{5.3a}
\bu_m=\bu_0,\qquad \bB_m=\bB_0,\quad  \mbox{in}\quad \R^2\times \{0\}.
\ee
Here  $(\bu_m,  p_m, \bB_m) \in  V^2_{\rm div} (Q)\times  L^{2}(Q) \times V^2_{\rm div}(Q)$ is called a weak solution to (\ref{5.2})-(\ref{5.1})
if  
\begin{align}
  & \intl_{Q} (
- \bu_m\cdot \partial _t\bfvarphi + \nabla \bu_m : \nabla \bfvarphi - 
\bu_m \otimes \bu_m : \nabla \bfvarphi  ) dx dt 
\cr
&\quad = \intl_{Q}  p_m \nabla \cdot  \bfvarphi  dx dt +
\intl_{Q} \Big((\nabla \times \bB_m) 
\times  \frac {\bB_m} {1+ \delta _m |\bB_m|} \Big) \cdot \bfvarphi  dxdt  + 
\intl_{Q} \bbf\cdot \bfvarphi  dxdt,   
\label{5.4}
\\
  &\intl_{Q} (- \bB_m \cdot \partial _t \bfvarphi + 
\nabla \bB_m : \nabla \bfvarphi ) dx dt
\cr
&\quad = -\intl_{Q} \Big((\nabla \times \bB_m- \bu_m) \times 
\frac {\bB_m} {1+ \delta _m |\bB_m|}\Big) \cdot \nabla \times \bfvarphi  dx  
+ \intl_{Q} \bg \cdot \nabla \times  \bfvarphi  dx  dt
\label{5.5}
\end{align} 
for all $\bfvarphi \in  C^\infty _{\rm c}(Q) $.   

\hspace*{0.5cm}
The existence of weak solutions to \eqref{5.2}--\eqref{5.3a} is given by the following 

\begin{lem}
\label{lem5.1}
Let $\bu_0 \in  L^2_{\rm div}, \bB_0 \in L^2 $ and  $\bbf, \bg \in  L^2(Q)$. Then for every $m\in \N$ there exists a weak solution 
$(\bu_m , p_m, \bB_m) \in V^2_{\rm div}(Q)\times L^2(0, T; L^2_{\rm loc} ) 
\times V^2_{\rm div}(Q) $  to  \eqref{5.2}--\eqref{5.3a}, such that 
\be\label{5.5a}
\nabla \bu_m \in  V^2(Q_r),\quad \forall\, \overline{Q} _r \subset Q. 
\ee 
Furthermore, this solution fulfills the energy equality
\begin{align}
 \frac {1} {2} \|\bu_m(t)\|_2^2 +  \frac {1} {2} \|\bB_m(t)\|_2^2 
+ \intl_{0}^{t} (\|\nabla \bu_m(s)\|_2 ^2 + \|\nabla \bB_m(s)\|_2 ^2) d s 
\cr
= \frac {1} {2} \|\bu_0\|_2^2 +  \frac {1} {2} \|\bB_0\|_2^2 + 
\intl_{0}^{t} \intl_{\R^2}  ( \bbf \cdot \bu_m + \bg \cdot \nabla \times \bB_m  ) dx ds  
\label{5.6}
\end{align}
for a.\,e. $t\in  (0,T)$. 

\end{lem}

{\em Proof }  Let $m\in \N$ be fixed.  Let $\beta _l \rightarrow 0^+ $ as $l \rightarrow + \infty $. 
By using the well-known monotone operator theory we get a weak solution $(\bu_{m, l}, p_{m, l}, \bB_{m, l} ) \in V^2_{\rm div}(Q)\times L^2(0,T; L^2_{\rm loc}) \times V^2_{\rm div}(Q) $ of the 
following   approximate system 
\begin{align} 
\partial _t \bu_{m,l}  +  &\frac {\bfomega_{m,l}} {1+\delta _m|\bB_{m,l}|+
 \beta _l |\bV_{m,l} |} \times \bu_{m,l} - \Delta \bu_{m,l} 
\cr
&\qquad \qquad = -\nabla p_{m,l} +
 (\nabla \times \bB_{m,l}) \times  \frac {\bB_{m,l}} {1+ \delta_m |\bB_{m,l}|+ 
\beta _l |\bV_{m,l} |} + \bbf,
\label{5.6b}
\\
\partial _t \bB_{m,l}+ &\nabla \times 
\frac {\bB_{m,l}} {1+\delta _m|\bB_{m,l}|+\beta _l |\bV_{m,l} |} \times \bu_{m,l}-  \Delta \bB_{m,l} 
\cr
&\qquad \qquad = 
-\nabla \times \Big (\nabla \times \bB_{m,l}\times  
\frac {\bB_{m,l} } {1+ \delta_m |\bB_{m,l}|+ \beta _l |\bV_{m,l} |} \Big)+  \nabla \times \bg
\label{5.6c}
\\
&\qquad \qquad\nabla \cdot \bu_{m, l}  =0,\quad \qquad\nabla \cdot \bB_{m, l}  =0\label{5.6a}  
\end{align}
in $Q=\R^2\times  (0,T)$  together with the initial condition 
\be\label{5.6d}
\bu_{m,l}=\bu_0,\qquad \bB_{m,l}=\bB_0,\quad  \mbox{in}\quad \R^2\times \{0\},
\ee
where
\[
\bV_{m, l} = \bfomega _{m, l} + \bB_{m, l}.    
\]
Clearly, the energy equality \eqref{5.6} holds true with 
$\bu_{m, l}$ in place of $ \bu _m$ 
and $\bB_{m, l} $ in place of  $\bB_{m}$ respectively.  
In particular, both $\{\bu_{m, l}\}$ and $\{\bB_{m, l}\}$ are bounded in $V^2(Q)$. 
Thus, by a standard reflexivity argument along with Banach-Alaoglu's compactness lemma,
eventually passing to a subsequence, we may assume there exist 
$\bu_m \in V^2_{\rm div} (Q)$ and $\bB_m \in   V^2_{\rm div}(Q)$  such that  
\begin{align}
 \nabla \bu_{m,l}  & \rightarrow \nabla \bu_m, 
\quad   \nabla \bB_{m, l}  \rightarrow \nabla \bB_m \quad \mbox{{\it weakly in}}\quad L^2(Q), 
\label{5.7}\\
\bu_{m,l}  &\rightarrow \bu_m, 
\quad   \bB_{m, l}  \rightarrow \bB_m \quad \mbox{{\it weakly$^\ast$ in}}\quad 
L^\infty (0,T; L^2)   \quad \mbox{as}\quad l \rightarrow +\infty.
\label{5.8}
\end{align}
Furthermore, by  Lions-Aubin's compactness lemma we see that 
\be\label{5.9}
\bu_{m,l}   \rightarrow \bu_m, 
\quad   \bB_{m,l}  \rightarrow \bB_m \quad \mbox{{\it strongly in}}\quad L^2(Q)
\quad\mbox{as}\quad l \rightarrow +\infty.  
\ee
Hence, thanks to \eqref{5.7}, \eqref{5.8} and \eqref{5.9} we are in a position to carry out 
the passage to the limit   
$l \rightarrow +\infty $ in the weak formulation of \eqref{5.6b}--\eqref{5.6a}. 
Accordingly, there exists  $p_m \in  L^2(0, T; L^2_{\rm loc} )$ such that 
$(\bu_m, p_m, \bB_m)$  is a weak solution to \eqref{5.2}--\eqref{5.3a}.   Verifying that $\bu_m$ and $\bB_m$ satisfying the energy equality \eqref{5.6},  it follows that   
\be\label{5.9b}
\nabla \bu_{m,l}   \rightarrow \nabla \bu_m, 
\quad   \nabla \bB_{m,l}  \rightarrow \nabla \bB_m \quad \mbox{{\it strongly in}}\quad L^2(Q)
\quad\mbox{as}\quad l \rightarrow +\infty.  
\ee
As $V^2(Q) \hookrightarrow L^4(Q)$ from \eqref{5.9b} we infer 
\be\label{5.9c}
\bu_{m,l}   \rightarrow \bu_m, 
\quad   \bB_{m,l}  \rightarrow \bB_m \quad \mbox{{\it strongly in}}\quad L^4(Q)
\quad\mbox{as}\quad l \rightarrow +\infty.  
\ee

\hspace*{0.5cm}
Next, applying $\nabla \times $ to both sides of \eqref{5.6b} and combining the result with  \eqref{5.6c}, we are led to 
\be\label{5.10}
\partial _t \bV_{m, l}  -  \Delta \bV_{m, l} = 
- \nabla \times  \Big(\frac {\bV_{m, l} } {1+ \delta_{m} |\bB_{m,l} |+ 
\beta _{l} |\bV_{m,l} |}\times  \bu_{m,l} \Big)  
+ \nabla \times \bh  \quad \mbox{in}\quad Q,
\ee
where $ \bh= \bg + \bbf $.  By using a routine smoothing argument one gets 
$\bV_{m, l} \in  V^2(Q_r) $ for all $\overline{Q} _r  \subset  Q$. 

\hspace*{0.5cm}
Now, let $\overline{Q} _r= \overline{Q_r(X_0)} \subset  Q$ be arbitrarily chosen. 
Let $\theta \in  C^\infty _{\rm c}(B_r \times (t_0-r^2, t_0]) $  be a test function suitable for 
$Q_{r/2} $.  Testing \eqref{5.10} by $\theta ^2\bV_{m, l} $,  we get  
\begin{align}
&\frac {1} {2} \intl_{B_r} \theta ^2(t)|\bV_{m, l}(t) |^2 dx +    \intl_{t_0- r^2}^{t} \intl_{B_r} 
\theta ^2 |\nabla \bV_{m, l} |^2 dx ds  
\cr    
&\quad  = \frac {1} {2} \intl_{t_0- r^2}^{t} \intl_{B_r} (\partial _t \theta ^2 + \Delta \theta ^2) |\bV_{m, l} |^2 dxds
\cr
 & \qquad-\intl_{t_0- r^2}^{t} \intl_{B_r}
\Big(\frac {\bV_{m, l} } {1+ \delta_{m} |\bB_{m,l} |+ 
\beta _{l} |\bV_{m,l} |}\times  \bu_{m,l} - \bh\Big)\cdot \nabla \times (\theta ^2\bV_{m,l}) dx ds 
\label{5.10a}
\end{align}
for a.\,e. $t\in  (t_0-r^2, t_0)$. 
 From the above identity using the embedding $V^2(Q_r) \hookrightarrow L^4(Q_r)$,
 it is readily seen that 
\begin{align}
&\bigg( \intl_{Q_r} \theta ^4 |\bV_{m,l}|^4  dxdt\bigg)^{1/2}  
\cr
& \le c\esssup_{t\in (t_0-r^2, t_0)} \intl_{B_r} \theta ^2(t) |\bV_{m,l} (t)|^2  dx 
 +   c\intl_{Q_r} \theta ^2|\nabla \bV_{m,l} |^2  + r^{-2}  |\bV_{m,l} |^2 + |\bh|^2 dx dt
\cr
& \le cr^{-2} (1+\|\bu_{m,l } \|^2_{4})  \intl_{Q_r} |\bV_{m,l} |^2 dx dt + 
c \|\bh\|_{2}^2 
\cr
& \hspace{3cm}  + 
{\hat C}\|\bu_{m,l} \|_{4, Q_r}   \bigg(\intl_{Q_r} \theta^4 |\bV_{m, l} |^4 dxdt\bigg)^{1/2},
\label{5.19}
\end{align}
with an absolute constant $ {\hat C}>0$. 
As $\bu_m \in L^4(Q)$, we may choose $0< r<  \sqrt{t_0}$ such that 
${\hat C} \|\bu_{m} \|_{4, Q_r} \le \frac {1} {4}$. Observing \eqref{5.9c}, there exists $l_0 \in \N$ 
such that ${\hat C} \|\bu_{m,l} \|_{4, Q_r} \le \frac {1} {2}$ for all $l\ge l_0$. Accordingly, \eqref{5.19} 
implies 
\begin{align}
&\bigg( \intl_{Q_r} \theta ^4 |\bV_{m,l}|^4  dxdt\bigg)^{1/2}  
\le cr^{-2} (1+\|\bu_{m,l } \|^2_{4})  \intl_{Q_r} |\bfomega _{m,l}+  \bV_{m,l} |^2 dx dt
+ c \|\bh\|_2^2
\label{5.20}
\end{align}
for $l\ge l_0$.  Since the right of \eqref{5.20} is bounded independently of $l\in \N$, by a 
constant $C(\bu_0, \bB_0, \bbf, \bg)$  by virtue of the lower semi continuity of the norm 
from \eqref{5.20} together with \eqref{5.10a} and  \eqref{5.19}  we get 
\be\label{5.21}
\|\nabla \bV_m\|_{2, Q_{r/2} } + \|\bV_m\|_{L^\infty (t_0-r^2/4, t_0; L^2(B_{ r/2}))}  + \|\bV_m\|_{4, Q_{r/2} } \le C(\bu_0, \bB_0, \bbf, \bg),
\ee
where $ \bV_m = \bB_m + \nabla \times \bu_m$.  By applying a 
standard covering argument, since $ \bB_m\in L^2$  we see that $ \nabla \times \bu _m \in V^2(Q_r)$ for all 
$ \overline{Q} _r \subset Q $. Whence, the assertion follows from the inequality 
\[
\| \nabla  \bu_m \|_{ 2, Q_{ r/2}} \le c r^{ -1} (\| \bu_m \|_{ 2, Q_r} + \| \nabla \times \bu _m\|_{ 2, Q_r})
\]
which completes the proof of the lemma. 
 \hfill \Beweisende

\vspace{0.3cm}
\hspace{0.5cm}
Next, we are going to carry out the passage to the limit $m \rightarrow +\infty $, which can be done by an analogous argument used in the proof of 
Lemma\,\ref{lem5.1}.  Observing the  energy equality \eqref{5.6}, we find that 
both $\{\bu_{m} \}$ and $\{\bB_{m}\} $ are bounded in $V^2(Q)$.  Eventually passing to a subsequence, we get the existence of 
 $\bu, \bB\in V^2_{\rm div} (Q)$  such that  
\begin{align}
 \nabla \bu_{m}  & \rightarrow \nabla \bu, 
\quad   \nabla \bB_{m}  \rightarrow \nabla \bB \quad \mbox{{\it weakly in}}\quad L^2(Q), 
\label{5.22}
\\
\bu_{m}  &\rightarrow \bu, 
\quad   \bB_{m, l}  \rightarrow \bB_m \quad \mbox{{\it weakly$^\ast$ in}}\quad 
L^\infty (0,T; L^2)   \quad \mbox{as}\quad m \rightarrow +\infty.
\label{5.23}
\end{align}
Furthermore, by  Lions-Aubin's compactness lemma we see that 
\be\label{5.24}
\bu_{m}   \rightarrow \bu, 
\quad   \bB_{m}  \rightarrow \bB \quad \mbox{{\it strongly in}}\quad L^2(Q)
\quad\mbox{as}\quad m \rightarrow +\infty.  
\ee
With the aid of \eqref{5.22}, \eqref{5.23} and \eqref{5.24} we are in a position  
to carry out the passage to the limit  
$m \rightarrow +\infty $  in the weak formulation of \eqref{5.2}--\eqref{5.3a}, which 
yields a weak solution $(\bu, p, \bB)$  to \eqref{1.2}--\eqref{1.3a}.   

\vspace{0.2cm}
\hspace{0.5cm}
Our next aim is to get a strong $L^4$ convergence of $\bu_m$. 

\begin{lem}
\label{lem5.2}
Let $\{(\bu_m, p_m, \bB_m)\}$ be  a sequence of weak solutions to \eqref{5.2}--\eqref{5.3a} obtained by  Lemma\,\ref{lem5.1}. 
Furthermore, suppose \eqref{5.22}--\eqref{5.24}. Then, for every $\overline{Q} _r \subset   Q$ 
there holds 
\be\label{5.25}
\bu_m \rightarrow \bu \quad \mbox{{\it strongly in}} \quad \bL^4(Q_r) \quad 
\mbox{as} \quad m \rightarrow +\infty.
\ee 
In addition, for every $X_0 \in Q$  there exists $0< r= r(X_0) < \sqrt{t_0}$ such that 
\begin{align}
&\|\nabla \bfomega _m\|_{2, Q_{r} } + \|\bfomega _m\|_{L^\infty (t_0-r^2, t_0; L^2(B_r))}  + 
\|\bfomega _m\|_{4, Q_{r} } 
\cr
&\qquad \le C(\bu_0, \bB_0, \bbf, \bg)\quad \forall\, m \in \N.  
\label{5.26}
\end{align}

\end{lem}

{\em Proof }
Let $ m\in \N$.  In view of Lemma\,\ref{lem5.1}, taking the sum of \eqref{5.2} and \eqref{5.3},   we see that $\bV_m = \bfomega _m + \bB_m \in  V^2_{\rm loc}(Q)$  is a weak solution to the following system   
\be\label{5.27}
\partial _t \bV_m - \Delta \bV_m = 
- \nabla \times  \Big(\frac {\bV_{m} } {1+ \delta_{m} |\bB_{m} |}
\times  \bu_{m} \Big)  + \nabla \times \bh \quad \mbox{in}\quad Q.
\ee
Here $ V^2_{\rm loc}(Q)$ contains all $  \bfvarphi \in L^2(Q)$ such that $
\varphi |_{ Q_r} \in V^2(Q_r)$ for all $ \overline{Q} _r \subset Q$.

Clearly, there exists $\bv_m \in  V^2_{\rm loc}(Q) $ such that $\nabla \times  \bv_m = \bV_m$. 
Thus,  from \eqref{5.27} we infer that 
\be\label{5.27a}
\partial _t \bv_m - \Delta \bv_m = - \nabla \pi _m
- \frac {\bV_{m} } {1+ \delta_{m} |\bB_{m} |}
\times  (\bv_m - \bb_m)+ \bh \quad \mbox{in}\quad Q,
\ee
where $\bb_m= \bv_m - \bu_m $.  By the definition of $\bv_m$ we have $\nabla \times \bb_m = \bB_m$.  

\hspace{0.5cm}
Let $\overline{Q} _r  \subset Q$ be fixed.  Eventually, replacing $\bv_m $ by $\bv_m(t)- (\bv_m(t))_{x_0, r} $ \,$(t\in t_0-r^2, t_0)$, observing \eqref{5.22}, \eqref{5.23} and \eqref{5.24}, by virtue of Sobolev's embedding theorem we easily verify that 
\begin{align}
\bV_m &\rightarrow \bV \quad \mbox{{\it weakly in}} \quad L^2(Q_r), 
\label{5.28}
\\
\bb_m &\rightarrow \bb \quad \mbox{{\it strongly in}} \quad L^6(Q_r) \quad 
\mbox{as} \quad m \rightarrow +\infty .
\label{5.30}
\end{align}
Indeed, we note that $|\bb_m(t)_{x_0, B_r}| = |\bu_m(t)_{x_0, B_r}| \le 
\|\bu_m\|_{L^\infty(0,T; L^2) }$. Consequently, by Sobolev-Poincar\'{e}'s inequality we see that 
$\|\bb_m\|_{q, Q_r} \le c \|\bu_m\|_{L^\infty(0,T; L^2) } + c \|\bB_m\|_{L^\infty(0,T; L^2) }$ 
for every $1\le q<+\infty $.
Once more appealing to \eqref{5.24}, eventually passing to a subsequence we may assume that 
\be\label{5.31}
\bB_m \rightarrow \bB \quad \mbox{{\it a.\,e. in}} \quad Q \quad 
\mbox{as} \quad m \rightarrow +\infty.
\ee
 By means of Vitali's convergence theorem, making use of \eqref{5.30} and \eqref{5.31}, we get 
\be\label{5.32}
\frac {\bb_m} {1+ \delta _m |\bB_m|} \rightarrow \bb \quad \mbox{{\it strongly in}} \quad L^6(Q_r) \quad 
\mbox{as} \quad m \rightarrow +\infty.
\ee
Next, we define the local pressure 
\begin{align*}
\nabla \pi _{m,1}&= \bE_{B_r} (\Delta \bv_m), 
\\
\nabla \pi _{m, 2} &= \bE_{B_r} \Big(-\frac {\bV_{m} } {1+ \delta_{m} |\bB_{m} |}
\times  (\bv_m - \bb_m)+ \bh\Big),
\\     
\nabla  \pi _{m,  \rm hm}  &= - \bE_{B_r}(\bv_m),  
\end{align*}
where $\bE_{B_r} : W^{ -1,\, q}(B_r) \rightarrow W^{ -1,\, q}(B_r) $ stands for the projection defined by the Stokes equation. Note that the restriction of $ \bE_{ B_r}$ to  $L^q(Q_r)$\, ($1<q<+\infty $) defines a projection in $ L^q(Q_r)$  
(cf. \cite{wol1, wol2} for details). We also note that $ \pi _{m,  \rm hm}(t)$ is harmonic 
in $ B_r$ for a.\,e.  all $ t\in (t_0-r^2, t_0)$. As it has been proved in \cite{wol1},\eqref{5.27a}  implies that the function
 $\bz_m= \bv_m + \nabla \pi_{m, \rm hm}   \in  V^2(Q_r)$ solves the following system in sense of distributions
\be\label{5.33}
\partial _t \bz_m - \Delta \bz_m = - \nabla (\pi_{m,1}+ \pi_{m,2})  
- \frac {\bV_{m} } {1+ \delta_{m} |\bB_{m} |}
\times  \bu_m + \bh\quad \mbox{in}\quad Q_r,
\ee
Let $\phi \in C^\infty_{\rm c} (Q_r)$ be a non-negative cut-off function. Testing \eqref{5.33} by $\phi  \bz_m$,
we obtain the following energy equality
\begin{align}
 &\intl_{Q_r} \phi  |\nabla \bz_m|^ 2 dxdt    
\cr
&\quad = \frac {1} {2}\intl_{Q_r} (\partial _t \phi + \Delta \phi ) |\bz_m|^2 dxdt 
+ \intl_{Q_r} \Big(\frac {\bV_{m} } {1+ \delta_{m} |\bB_{m} |}
\times  \bb_m + \bh \Big)\cdot \phi \bz_m dx dt    
\cr
&\qquad\qquad+ \intl_{Q_r} (\pi_{m, 1} + \pi _{m, 2}) 
\nabla \phi \cdot \bz_m dx dt.
\label{5.34}
\end{align}
Verifying 
\[
\Big\|\frac {\bV_{m} } {1+ \delta_{m} |\bB_{m} |}
\times  \bu_m\Big\|_{L^{3/2}(0, T; L^{6/5} ) } \le \|\bV_m\|_2 \|\bu_m\|_{L^6(0,T; L^3)}  
\le C(\bu_0, \ldots),
\] 
we may estimate the pressure $\pi _{m, 2} $ in $L^{3/2}(Q_r) $ by using  the 
Sobolev-Poincar\'{e} inequality as follows
\begin{align}
\|\pi _{m, 2} \|_{3/2, Q_r} &\le  c \|\nabla \pi _{m, 2} \|_{L^{3/2}(t_0-r^2, t_0; L^{6/5}(B_r) )}
\cr       
&\le c \Big\|\frac {\bV_{m} } {1+ \delta_{m} |\bB_{m} |}
\times  \bu_m\Big\|_{L^{3/2}(0, T; L^{6/5} ) } + c \|\bh\|_2 \le C(\bu_0, \ldots).  
\label{5.34a}
\end{align}
Furthermore, we immediately get 
\be\label{5.34b}
\|\pi _{m,1} \|_{2, Q_r} \le c \|\nabla \bv_m\|_{2} \le c \|\nabla \bu_m\|_2 +
c \| \bB_m\|_2  \le C(\bu_0, \ldots). 
\ee
Observing \eqref{5.24} along with \eqref{5.30}, we find 
\be\label{5.35}
\bv_{m }  \rightarrow \bv\quad \mbox{{\it strongly in}} \quad L^3(Q_r) 
\quad \mbox{as} \quad m \rightarrow +\infty,
\ee
where $\bv = \bu + \bb$.  Thus, having
\begin{align}
\nabla  \pi_{m, \rm hm} 
\rightarrow \nabla  \pi_{ \rm hm}   
\quad \mbox{strongly in }\quad L^3(Q_r)\quad  \mbox{as}\quad 
m \rightarrow \infty,
\label{5.35a}
\end{align}
 where $\nabla  \pi_{\rm  hm} = -\bE_{B_r}(\bv) $,   it follows that 
\be\label{5.36}
\bz_m\rightarrow \bz\quad \mbox{{\it strongly in}} \quad L^3(Q_r) 
\quad \mbox{as} \quad m \rightarrow +\infty.
\ee
Now, with help of \eqref{5.34a}, \eqref{5.34b} and \eqref{5.36} we get 
\[
\lim_{m \to \infty }\intl_{Q_r} (\pi_{m, 1} + \pi _{m, 2}) 
\nabla \phi \cdot \bz_m dx dt    = \intl_{Q_r} (\pi_{ 1} + \pi _{ 2}) 
\nabla \phi \cdot \bz dx dt, 
\]
where 
\begin{align*}
\nabla \pi _1 &= \bE_{B_r} (\Delta \bv), 
\\
\nabla \pi _2 &= \bE_{B_r} (-\bV
\times  \bu + \bh).
\end{align*}

\hspace{0.5cm}
On the other hand, making use of \eqref{5.32}, together with \eqref{5.24} and \eqref{5.36} 
we see that 
\[
\lim_{m \to \infty } \intl_{Q_r} \Big(\frac {\bV_{m} } {1+ \delta_{m} |\bB_{m} |}
\times  \bb_m + \bh \Big)\cdot \phi \bz_m dx dt    = \intl_{Q_r} 
(\bV \times  \bb + \bh )\cdot \phi \bz dx dt.   
\]
Furthermore, thanks to \eqref{5.36} we obtain 
\[
\lim_{m \to \infty }\frac {1} {2}\intl_{Q_r} (\partial _t \phi + \Delta \phi ) |\bz_m|^2 dxdt
= \frac {1} {2}\intl_{Q_r} (\partial _t \phi + \Delta \phi ) |\bz|^2 dxdt.
\]
Hence, we are in the position to carry out the passage to the limit 
 $ m \rightarrow +\infty $ in \eqref{5.34} to get 
\begin{align}
&\lim_{m \to \infty } \intl_{Q_r} \phi  |\nabla \bz_m|^ 2 dxdt   
\cr 
& \quad = \frac {1} {2}\intl_{Q_r} (\partial _t \phi + \Delta \phi ) |\bz|^2 dxdt 
+ \intl_{Q_r} (\bV \times  \bb + \bh)\cdot \phi \bz dx dt +
\intl_{Q_r} (\pi_{ 1} + \pi _{ 2}) 
\nabla \phi \cdot \bz dx dt.   
\label{5.38}
\end{align}
Accordingly, we see that $\bz\in  V^2(Q_r)$ and 
\be\label{5.40}
\partial _t \bz- \Delta \bz = - \nabla (\pi_{1}+ \pi_{2})  
- \bV\times  \bu + \bh \quad \mbox{in}\quad Q_r,
\ee
in sense of distributions.  Taking into account that $\bz \in L^4(Q_r)$, and $\bV\times  \bu  \in  L^{4/3}(Q) $, 
we obtain  the following energy equality
\begin{align}
 &\intl_{Q_r} \phi  |\nabla \bz|^ 2 dxdt    
\cr
&\quad = \frac {1} {2}\intl_{Q_r} (\partial _t \phi + \Delta \phi ) |\bz|^2 dxdt 
+ \intl_{Q_r} (\bV\times  \bb + \bh)\cdot \phi \bz dx dt    
\cr
&\qquad\qquad+ \intl_{Q_r} (\pi_{1} + \pi _{2}) 
\nabla \phi \cdot \bz dx dt.    
\label{5.40a}
\end{align}
Thus, observing \eqref{5.24},  combining \eqref{5.38} and \eqref{5.40a} using a well-known liminf-limsup argument noticing that 
$ \sqrt[]{\phi }\nabla \bz_m \rightarrow\sqrt[]{\phi }\nabla \bz$ 
weakly in $ L^2(Q_r)$ 
we  get 
\[
\sqrt[]{\phi } \nabla \bz_m \rightarrow \sqrt[]{\phi }\nabla \bz   
\quad \mbox{strongly  in }\quad L^2(Q_r)\quad \mbox{as}\quad m \rightarrow \infty.
\]
On the other hand, since $ \pi _{ \rm hm}$ is harmonic, thanks to \eqref{5.35a} we get
\begin{align*}
\sqrt[]{\phi }\nabla ^2 \pi_{m, \rm hm} 
\rightarrow \sqrt[]{\phi }\nabla ^2 \pi_{ \rm hm}   
\quad \mbox{strongly in }\quad L^2(Q_r)\quad  \mbox{as}\quad 
m \rightarrow \infty.
\end{align*}
As $ \nabla \bv_m = \nabla \bz_m - \nabla ^2 \pi_{m, \rm hm}$ a.\,e. in $ Q_r$, we arrive at
\[
\sqrt[]{\phi } \nabla \bv_m \rightarrow \sqrt[]{\phi }\nabla \bv   
\quad \mbox{strongly  in }\quad L^2(Q_r)\quad\mbox{as}\quad m \rightarrow \infty.
\]
Hence,   thanks to the embedding  $V^2(Q_r) \hookrightarrow L^4(Q_r)$ along with \eqref{5.30}  we get  
\[
\sqrt[]{\phi }\bu_m \rightarrow \sqrt[]{\phi}\bu   
\quad \mbox{strongly  in }\quad L^4(Q_r)\quad \quad \mbox{as}\quad m \rightarrow \infty.
\]
Since the above statement holds for any cylinder $Q_r \subset Q$,  we get the first claim \eqref{5.25} 
of the lemma.   

\hspace*{0.5cm}
Now, it remains to verify \eqref{5.26}. In fact, according to Lemma\,\ref{lem5.1}, we have $\bV_m \in  L^4(Q_r)$, which 
implies that $\frac {\bV_m} {1+ \delta _m |\bB_m|} \times  \bu_m \in  L^2(Q_r)$. This allows us to test 
\eqref{5.27} with $\theta ^2 \bV_m$, where $\theta \in  C^\infty _{\rm c} (B_r \times (t_0-r^2, t_0])$.  
Arguing as in the proof of Lemma\,\ref{lem5.1},  we obtain 
\begin{align*}
&\frac {1} {2} \intl_{B_r} \theta ^2(t)|\bV_{m}(t) |^2 dx +    \intl_{t_0- r^2}^{t} 
\intl_{B_r} \theta ^2 |\nabla \bV_{m} |^2 dx ds  
\\     
&\quad  = \frac {1} {2} \intl_{t_0- r^2}^{t} \intl_{B_r} (\partial _t \theta ^2 + \Delta \theta ^2) |\bV_{m} |^2 dxds
\\
 & \qquad-\intl_{t_0- r^2}^{t} \intl_{B_r}
\Big(\frac {\bV_{m} } {1+ \delta_{m} |\bB_{m} |}\times  \bu_{m} - \bh\Big)\cdot 
\nabla \times (\theta ^2\bV_{m}) dx ds 
\end{align*}
for a.\,e. $t\in  (t_0-r^2, t_0)$, which leads to 
\begin{align}
&\bigg( \intl_{Q_r} \theta ^4 |\bV_{m}|^4  dxdt\bigg)^{1/2}  
\cr
&\quad \le cr^{-2} (1+\|\bu_{m} \|^2_{4})  C(\bu_0, \ldots) + 
{\hat C}  \|\bu_{m} \|_{4, Q_r}   \bigg(\intl_{Q_r} \theta^4 |\bV_{m} |^4 dxdt\bigg)^{1/2}. 
\label{5.41}
\end{align}
Whence, the proof  of \eqref{5.26} can be  completed  by a similar argument to the proof of 
Lemma\,\ref{lem5.1}, by using the strong $L^4$ convergence \eqref{5.25}. \hfill \Beweisende

%
%
\section{Proof of Theorem\,\ref{thm1.1}} 
\setcounter{secnum}{\value {section}
\setcounter{equation}{0}
\renewcommand{\theequation}{\mbox{\arabic{secnum}.\arabic{equation}}}}

Let $(\bu_m, p_m, \bB_m) \in V^2_{\rm div}(Q)\times  L^2(0,T; L^2_{\rm loc}) \times  V^2_{\rm div}(Q)$ 
be a weak solution to the approximate system \eqref{5.2}--\eqref{5.3a}  such that  
$\nabla \bu_m \in V^2_{\rm loc}(Q) $ $ (m\in \N)$, 
which can be guaranteed by Lemma\,\ref{lem5.1} (for the definition of $ V^2_{\rm loc}(Q)$ see Section\,4). 

\hspace*{0.5cm}
In our discussion below we use the following notation. 
Let  $X_0 = (x_0, t_0)\in  Q$. 
\begin{align*}
E_m (r ) =E_m (r, X_0) &:= \bigg(\intmw_{Q_r(X_0)} |\bB_m- (\bB_m)_{r, X_0} |^4 dx dt\bigg)^{1/4},   
\\     
F_m (r ) =F_m (r, X_0) & := \bigg(r^{-2} \intl_{Q_r(X_0)} |\nabla \bB_m|^2 dx dt\bigg)^{1/2},
\\[0.2cm]
H_m (r ) =H_m (r, X_0) & := \|\bu_m\|_{4, Q_r} + r^{-1} \|\bg\|_{2, Q_r} 
\quad   0< r< \sqrt[]{t_0}.
\end{align*}

\hspace{0.5cm}
Next, we  define  the set of possible singularities of $ B$ 
by means of  $ \Sigma(\bB)= \cup _{ k=1}^\infty \Sigma_k \cup \Sigma _\infty  $, where 
\begin{align*}
\Sigma_k &:= \bigcup _{ 0< \rho <T} \bigcap _{ 0<r \le  \rho }  
\bigg\{X_0 \in  \R^2\times (r , T) \bigg|  \liminf_{m \to \infty} 
F_m(r, X_0) \ge \frac {1} {k} \bigg \},\quad  k\in \N,
\\
\Sigma _\infty &:=\Big\{ X_0 \in Q \,\Big|\, \sup_{0<r<\sqrt{t_0}}  
|\bB_{r, X_0}|=+\infty  \Big\}.
\end{align*}

\hspace*{0.5cm}
Let $ Q_r =Q_r(X_0)\subset Q$ be any  cylinder such that  condition \eqref{3.11} is fulfilled for 
$\bB= \bB_m$ and $\bu=\bu_m$, i.\,e.  
\be\label{6.11}
C_1  \Big\{F_m(r)^2 + \|\bu_m\|^2_{4, Q_r}\Big\} \le \frac {1} {2}.    
\ee 
As stated in Remark\,\ref{rem3.4}, the condition \eqref{3.11} implies \eqref{3.14}. Thus,  
\eqref{6.11} implies 
\begin{align}
 E_m(r/2) &\le  
C_2  (1+ |(\bB_m)_{r, X_0} |^2 ) \Big\{F_m(r) + F_m(r)^2 + H_m(r) + H_m(r)^2\Big\}.
\label{6.12}
\end{align}
On the other hand, \eqref{3.15}  with $ \bB =\bB_m$  and 
$\bu= \bu _m$  
reads
\begin{align}
F_m(r/2) &\le  
C_3   (1+ |(\bB_m)_{r, X_0} | ) \Big\{E_m(r) + E_m(r)^2 + H_m(r) + H_m(r)^2\Big\}.
\label{6.10}
\end{align}

\hspace*{0.5cm}
Let $X_0 \in  Q\setminus  \Sigma (\bB)$  be fixed.  Set $d_0 = \sqrt[]{t_0}/2$ 
and $K = \overline{Q_{d_0}}$. 
Appealing to Lemma\,\ref{lem5.2},  and applying Sobolev's embedding theorem, 
we see that  
\be\label{6.10a}
\|\bu_m\|_{ 8/(4- \lambda ), K} \le L\quad \forall\, m \in \N,    
\ee
where  $L=\const>0$ depends on $d_0, \bu_0, \bB_0, \bbf$ and $\bg$ only.  
Furthermore,  we may choose $0<R_1< d_0$ such that 
\be\label{6.10b}
C_1 \|\bu_m\|_{4, Q_{R_1} }^2 \le \frac {1} {16}\quad \forall\, m\in \N,   
\ee
where $C_1$ stands for the constant appearing in \eqref{6.11}.  
Using H\"older's inequality, recalling the assumption on $\bg$  along with 
\eqref{6.10a},  it follows that 
\begin{align}
H_m(r,X_0) & \le  \Big(\pi ^{\lambda /8-1/4} \|\bu_m\|_{8/(4-\lambda ), K} +  
\|\bg\|_{{\cal M}^{2, \lambda }(K) }\Big)  r^{(\lambda-2) /2}
\cr     
&\le C_4  r^{(\lambda -2)/2} \quad \forall\, \quad 0<r\le R_1.  
\label{6.10c}
\end{align}

\hspace*{0.5cm}
Next, we set 
\[
M:= 512\sup_{0< r< d(X_0/2)} (|\bB|)_{r, X_0 }+1 <+\infty. 
\]
Let $ 0< \alpha < \frac{2-\lambda }{2}$. We take  $\tau  > 0$  such that 
\be\label{6.18}
2 \tau^{1-\alpha }   C_0 (1+ M^5) \le \frac {1} {2} \quad \mbox{and}\quad 
\tau ^{\alpha } \le \frac {1} {2}
\ee
(Recall, the constant $ C_0>0$ has been defined in Lemma\,\ref{lem4.1}).

\hspace{0.5cm}
Now, let $\var _0=\var_0(\tau , M, L,  K, \alpha ),  R_0=R_0(\tau , M, L, K, \alpha )$ 
and $\delta _0=\delta _0(\tau , M, L, K, \alpha )$ denote the numbers according to 
Lemma\,\ref{lem4.2}.  
In addition, we define $\var _1>0$ by the relation 
\be\label{6.20b}
2 \tau ^{-4} \var _1 = 1. 
\end{equation}

Next we may choose  $0 < R_2 \le \min\{R_0, R_1\}$ such that  the following conditions hold
\begin{align}
& C_2(1+ M^2)(C_4+C_4^2) R_2^{(\lambda -2)/2} \le 
\frac {1} {8}\min\{\var _0, \var _1\},
\label{6.30a}
\\
& \qquad \qquad  2R^{\alpha }_2 \le   \frac {1} {2} \min\{\var _0, \var _1\}.
\label{6.23}
\end{align} 

Now, we take  $k\in \N$  such that 
\begin{align}
C_2 (1+ M^2) \Big\{\frac {1} {k}  + \frac {1} {k^2}\Big\}
\le \frac {1} {8} \min\{\var _0, \var _1\} \quad \mbox{and}\quad \frac {C_1} {k}\le \frac {1} {4}.
\label{6.20a}
\end{align}

Owing to  $ X_0\in Q \setminus \Sigma _k$ eventually replacing  $R_2$ by a smaller number we may also assume that  $ \liminf_{m \to \infty} F_m (R_2, X_0) < \frac1k$. 
Accordingly we are able to  select 
a subsequence $ \{ m_j\}$ such that 
\be\label{6.23a}
F_{m_j} (R_2, X_0 ) <  \frac {1} {k}\quad \forall\,j \in \N.
\ee

Since $\bB_m \rightarrow \bB$ in $L^1(Q_{R_2} )$ 
and $ \delta _m \rightarrow  0$ as $ m \rightarrow +\infty$, there exists $m_0\in \N $ with the property  
\be\label{6.19}
(|\bB_{m}|) _{R_2 , X_0} \le (|\bB|) _{R_2 , X_0} + \frac{1}{512} 
\le
\frac {M} {512} \quad  
 \mbox{and}\quad  \delta _{m} \le \delta _0 \quad \forall\, m \ge m_0.  
\ee

\hspace{0.5cm}
Observing \eqref{6.23a}, \eqref{6.20a} and \eqref{6.10b},  we have 
\be\label{6.10d}
C_1 \Big\{F_{m_j} (R_2  , X_0) +4 \|\bu_{m_j} \|_{4, Q_{R_2 }(X_0) }\Big\} \le \frac {1} {2}\quad  \forall\, j\in \N.
\ee
As \eqref{6.10d}  implies \eqref{6.12},  employing  
\eqref{6.19}, \eqref{6.23a} and  \eqref{6.10c},  we get  
\begin{align}
 E_{m_j} (R_2 /2, X_0)  &\le C_2 (1+ M^2) \Big\{\frac {1} {k}  + \frac {1} {k^2}
+ (C_4+C_4^2) R_2^{(\lambda -2)/2} \Big\} 
\label{6.20}
\end{align}
for all $m_j \ge m_0$.   
In view of \eqref{6.20a} and \eqref{6.30a}, \eqref{6.20} 
gives 
\begin{align}
 E_{m_j} (R_2 , X_0)  &\le \frac {1} {4} \min\{\var _0, \var _1\}\quad \forall\, m_j \ge m_0. 
\label{6.20c}
\end{align}

\hspace{0.5cm}
Set $ R_3= R_2 /2$. Let $Y \in  Q_{R_2}(X_0) $.  Clearly, 
\begin{align}
E_{m_j}(R_3, Y)  &\le 2 E_{m_j}( R_2 ,  X_0) \le \frac {1} {2} \min\{\var _0, \var _1\}, 
\label{6.21}
\\
|(\bB_{m_j})_{R_3, Y}| &
\le  256 \intmw_{Q_{R_2}(X_0) } |\bB_{m_j}|dxdt \le \frac {M} {2}.
\label{6.22}
\end{align}

We claim that for every $i \in  \N\cup \{0\}$,  there holds 
\begin{align}
E_{m_j} (\tau ^{i} R_3, Y) &\le 2^{- i} \tau^{\alpha i} E_{m_j} (R_3, Y) + (1- 2^{-i} ) \tau ^{\alpha i} R_3^{\alpha },
\label{6.24}  
\\
|(\bB_{m_j}) _{ \tau ^i R_3, Y}|&\le  M - 2^{-i+1}.
\label{6.25}
\end{align}
In fact, for $i=0$, \eqref{6.24} is trivially fulfilled,  while \eqref{6.25} holds in view of \eqref{6.22}.  

\hspace*{0.5cm}
Now, we assume that both  \eqref{6.24} and \eqref{6.25} are fulfilled for $i\in  N\cup \{0\}$.  
Then \eqref{6.24} together  with \eqref{6.21} and \eqref{6.23}   implies
\be\label{6.26}
E_{m_j} (\tau ^{i} R_3, Y) + \tau^{\alpha i } R_3^{\alpha }  \le \tau^{\alpha i} 
 (E_{m_j} (R_3, Y) + 2R_3^\alpha)
\le \tau^{\alpha i} \min\{\var _0, \var _1\}.
\ee 
In particular,  observing \eqref{6.25} we have  
\[
E_{m_j} (\tau ^{i} R_3, Y) + (\tau^{ i} R_3)^{\alpha}\le \var _0,\quad 
|(\bB_{m_j}) _{\tau ^i R_3, Y}|\le  M. 
\]
Thus, we are in a position to apply Lemma\,\ref{lem4.2} with $R= \tau ^i R_3$. This together with 
\eqref{6.24}  gives 
\begin{align}
E_{m_j} (\tau^{i+1} R_3, Y) &\le 2 \tau C_0 (1+ M^5) (E_{m_j} (\tau^{i} R_3, Y) + 
\tau ^{\alpha i} R_3^{\alpha }) 
\cr
&\le \frac {1} {2} \tau^{\alpha} E_{m_j} (\tau^{i} R_3, Y) + \frac {1} {2}\tau ^{\alpha (i+1)} R_3^{\alpha }
\cr
&\le 2^{-(i+1)} \tau^{\alpha (i+1) } E_{m_j} (R_3, Y) + (1- 2^{-(i+1)} ) \tau ^{\alpha (i+1)} R_3^{\alpha }.
\label{6.27}
\end{align}
Consequently \eqref{6.24} holds true for $i +1$.  

\hspace{0.5cm}
Now, it remains to show \eqref{6.25} for $i+1$.  First, from 
\eqref{6.24} along with  \eqref{6.21} and \eqref{6.23} we infer 
\be\label{6.28}
E_{m_j} (\tau ^{i} R_3, Y) \le \tau^{\alpha i} (E_{m_j} (R_3, Y) + R_3^{\alpha }) 
\le \tau ^{\alpha i} \var _1.  
\ee
Using the triangle inequality and Jensen's inequality,  we find 
\begin{align*}
  |(\bB_{m_j}) _{\tau^{i+1} R_3, Y}| &\le |(\bB_{m_j}) _{ \tau^{i}  R_3, Y}| + 
\Big|(\bB_{m_j}) _{\tau^{i+1}  R_3,  Y}- (\bB_{m_j}) _{\tau^{i}  R_3, Y}\Big|
\\     
&\le |(\bB_{m_j}) _{\tau^{i}  R_3, Y}| + 2 \tau ^{-4} E_{m_j} ( \tau ^{i} R_3, Y).  
\end{align*}
Estimating the first term on the right by using \eqref{6.25} and the second one 
by the aid of \eqref{6.28}  together with \eqref{6.18} and \eqref{6.20b}, we obtain 
\begin{align*}
|(\bB_{m_j}) _{ \tau^{i+1}  R_3, Y}|  & \le M- 2^{-i+1} + 2 \tau ^{-4} \tau ^{\alpha i} \var _1
\\     
&\le M-2^{-i+1} + 2^{-i}= M- 2^{-i}.  
\end{align*}
This completes the proof of \eqref{6.25} for $i+1$. Whence, the claim.

\hspace*{0.5cm}
Since \eqref{6.24} holds true for every $ Y\in \overline{Q_{R_3}(X_0)}$, 
by a standard iteration argument we get a constant $C_5>0$ such that 
\be\label{6.31a}
\bigg(\intmw_{Q_r(Y)} |\bB_{m_j} - (\bB_{m_j}) _{r, Y} |^4  dx\bigg)^{1/4} \le  C_5 r^{\alpha } \quad \forall\, 0<r < R_3, \quad \forall\, Y\in  \overline{Q _{R_3}  (X_0)}.
\ee
Thus, by means of the lower semi continuity of the $L^4$-norm the above inequality remains true for $\bB$.  Using the well-known integral characterization of the 
 H\"{o}lder continuity  in the parabolic setting\cite{dap}, we obtain 
\be\label{6.31}
\bB|_{\overline{Q_{R_3}  (X_0)}} \in  C^{\alpha, \alpha /2 }(\overline{Q_{R_3} (X_0)})
\ee
(For the definition of $C^{\alpha, \alpha /2 }(\overline{Q_{R_3} (X_0)}) $ see appendix below).
Clearly, \eqref{6.31a} shows that 
\[
\lim_{r \to 0^+} E_{m_j} (r, Y) = 0\quad \mbox{{\it uniformly for }}\quad Y\in Q_{R_3}(X_0) \quad \text{and} \quad  j\in \N.  
\]
Hence, in view of \eqref{6.10} we get $
Y \not\in  \bigcup_{k=1}^\infty  \Sigma _{k}. $
Taking into account that $\bB$ is H\"older continuous on 
$\overline{Q_{\rho _0}(X_0)} $, it follows that $ Q_{\rho _0}(X_0) 
\subset Q \setminus \Sigma _\infty $ and thus
\[
Q_{R_3}(X_0) \subset Q \setminus \Sigma (\bB).
\]
Consequently, $\Sigma (\bB)$ is a closed set. 
This completes the proof of the main theorem. \hfill \Beweisende

\begin{thm}
\label{thm6.1}
For the singular set constructed in the proof of Theorem\,\ref{thm1.1} we have 
\be\label{6.32}
d{\cal P}_{\beta }(\Sigma (\bB))=0 \quad \forall\, \beta >2,
\ee
where $d{\cal P}_{\beta }(\cdot )$ is the $\beta-$dimensional parabolic Hausdorff measure. In particular, the Hausdorff dimension of $ \Sigma(\bB)$ satisfies 
$\dim_{{\cal H}} (\Sigma (\bB))\le 2$. 
\end{thm}

{\it Proof } Let $2<\beta \le  \lambda $ be arbitrarily chosen.  First we show that 
\[
d {\cal P}_\beta (\Sigma _k)=0 \quad \forall\, k\in \N .
\]

\hspace{0.5cm}
Let $X_0 \in  \Sigma _k$. Fix $\var >0$.  Then there exists 
$0<r (X_0) <\varepsilon $ and $m(X_0)\in  \N$, such that 
\be\label{6.33}
r(X_0)^{-2}  \intl_{Q_{r(X_0)}(X_0) }  |\nabla \bB_m|^2 dxdt  \ge \frac {1} {2k} \quad \forall\, m\ge m(X_0).
\ee
Clearly, the family of cylinders $\{Q_{r(X_0)}(X_0) \}_{X_0 \in \Sigma _k} $ forms  a covering of 
$\Sigma _k$.  Thanks to the Vitali covering lemma there exists a pairwise disjoint 
family $\{ Q_{r_i}(X_i)\}_{i\in \N} $  $(r_i := r(X_i)) $ such that 
$ \{Q_{3 r_i}(X_i)\}_{i\in \N}$  covers  $ \Sigma _k$. Let $N\in \N$ be arbitrarily chosen. Set 
\[
m_N : = \max \{m(X_1), \ldots, m(X_N)\}.  
\] 
Then, from \eqref{6.33} with $X_0 =X_i$\, ($i=1,\ldots,N$) and  $m=m_N$ we infer  
\begin{align*}
\suml_{i=1}^{N}  r_i^{\beta}  &\le
\var ^{\beta -2} \suml_{i=1}^{N}  r_i^2 \le 2 \var ^{\beta -2} k  \suml_{i=1}^{N} \intl_{Q_{r_i}(X_i) }  |\nabla \bB_{m_N} |^2 dxdt \le 2
\var ^{\beta -2} k \intl_{Q}  |\nabla \bB_{m_N} |^2
\\     
&\le \var ^{\beta -2} k C(\|\bu_0\|_2, \ldots ). 
\end{align*}
This shows that 
\be\label{6.34}
\suml_{i=1}^{\infty }  r_i^{\beta }\le   \var ^{\beta -2} k C(\|\bu_0\|_2, \ldots ). 
\ee
Consequently, $d {\cal P}_{\beta }  (\Sigma _k)=0$, which implies that 
$d {\cal P}_{\beta }\Big(\bigcup_{ k=1}^\infty \Sigma_k\Big)=0 $.

\hspace*{0.5cm}
Now, it remains to prove that 
$d {\cal P}_{\beta } (\Sigma_\infty)=0$.  As we will see below this follows easily 
from the following implication 
\be\label{6.36a}
\sup_{0<r< \sqrt{t_0}} r^{-\beta } \intl_{Q_r(X_0)}  |\nabla B|^2 dxdt <+\infty      
\quad \Longrightarrow \quad X_0 \notin \Sigma _\infty,\quad  X_0\in Q. 
\ee
Indeed, let   $X_0 \in  Q$ such that  the  condition on the left in \eqref{6.36a} holds true.  
Choose $0<\rho _0< \sqrt{t_0}$ sufficiently small (specified below) and set $r_i = 2^{-i} \rho _0$ \, $(i\in \N)$. 

\hspace{0.5cm}
Fix $i\in \N$. By using the parabolic Poincar\'{e}-type inequality (see Lemma A.1, appendix below),  arguing as  in the proof of 
\eqref{3.10},  we estimate 
\begin{align}
& \intmw_{Q_{r_i} } |\bB- \bB_{r_i , X_0} |^2 dxdt 
\cr
&\quad \le c(1+ |\bB_{r_i, X_0} |^2) r_i^{-2} \intl_{Q_{r_i} } |\nabla \bB|^2  dx dt      
\cr
& \qquad  + c (1+  |\bB_{r_i, X_0} |^2)  r^{-2}\intl_{Q_{r_i} } (|\bg|^2 +|\bu|^2) dx dt 
\cr
 & \qquad + C_6 \Bigg\{r^{-2}_{i}   \intl_{Q_{r_i} } |\nabla \bB|^2 + 
\bigg(\intl_{Q_{r_i} }  |\bu|^4 dxdt \bigg)^{1/2}\Bigg\}  
\intmw_{Q_{r_i} }   |\bB- \bB_{r_i, X_0} |^2  dxdt   
\label{6.40}
\end{align}
for an absolute constant $C_6>0$.  Due to $\bu \in L^4(Q)$ and our assumption on $X_0$ we may choose $\rho _0$  sufficiently small such that the numerical value in $\{\ldots\} $ is less than $\frac {1} {2C_6}$, which leads to 
\begin{align}
& \intmw_{Q_{r_i} } |\bB- \bB_{r_i , X_0} |^2 dxdt 
\cr
&\quad \le 2c(1+ |\bB_{r_i, X_0} |^2) r_i^{-2} \intl_{Q_{r_i} } |\nabla \bB|^2  dx dt      
\cr
& \qquad \qquad  + 2c (1+  |\bB_{r_i, X_0} |^2)  r^{-2}_i\intl_{Q_{r_i} } (|\bg|^2 +|\bu|^2) dx dt.
\label{6.41}
\end{align}
Appealing to Lemma\,\ref{lem5.2},  we see that $\bu \in L^{q}_{\rm loc}(Q)  $ for all $1\le q <+\infty $. 
In particular, $\bu \in {\mathcal M}^{2, \lambda }(Q_{\sqrt{t_0}/2})$.  Recalling  that $\bg \in {\cal M}^{2, \lambda }(Q) $ 
and $\beta \le \lambda $ from \eqref{6.41},  we deduce that  
\begin{align}
& \intmw_{Q_{r_i} } |\bB- \bB_{r_i , X_0} |^2 dxdt \le c(1+ |\bB_{r_i, X_0} |^2)  r_{i}^{\beta -2}    
\label{6.42}
\end{align}
with a constant $c>0$ depending neither on $r_i$ nor on $\rho _0$.  
Using the triangle inequality and employing  \eqref{6.42},  it follows  that 
\be\label{6.43}
\Big||\bB_{r_{i+1} , X_0} | - |\bB_{r_i, X_0} |\Big|\le 
C_7 (1+ |\bB_{r_i, X_0} |)  r_{i}^{(\beta -2)/2},
\ee
where $C_7= \const >0$ is independent on $r_i$ and $\rho _0$.  
Thus,  eventually replacing $\rho _0$ by a smaller one,  we may assume that 
\[
C_7\suml_{i=0}^{\infty } r_i^{(\beta -2)/2} = C_7 \rho _0^{(\beta -2)/2} \frac {1} {1- 2^{(\beta -2)/2} } 
\le \frac {1} {2}.     
\]  
Then, with help of Lemma\,\ref{lemA.2} (see appendix below) from \eqref{6.43} we conclude that 
\be\label{6.44}
|\bB_{r_{i} , X_0} |\le 1+ 2 |\bB_{\rho _0 , X_0} | \quad \forall\, i\in \N,
\ee
what completes the proof of \eqref{6.36a}. 

\hspace*{0.5cm}
Now, let $\var >0$ be arbitrarily chosen. According to \eqref{6.36a} for every $X_0\in \Sigma _\infty $ 
we may choose $0<r=r(X_0) \le \var $ such that 
\[
r^{-\beta } \intl_{Q_r(X_0)}  |\nabla B|^2 dxdt \ge \frac {1} {\var }.
\]
Thus,  by the Vitali covering lemma there exists a pairwise disjoint family  $\{Q_{r_i}(X_i) \}$ $(r_i := r(X_i)$
such that $ \{Q_{3 r_i}(X_i)\} $ covers $\Sigma _\infty $. Similarly to the above we conclude
\[
\suml_{i=1}^{\infty } r_i^{\beta } \le c \var \|\nabla \bB\|^2_{2}.
\]
Thus, $d {\cal P}_\beta  (\Sigma _\infty )=0$,  and the proof of the theorem is complete. \hfill \Beweisende

\vspace{0.5cm}
$$
\mbox{\bf Acknowledgements}
$$
Chae was partially supported by NRF grants 2006-0093854 and  2009-0083521, while Wolf has been supported by the Brain Pool Project of the Korea Federation of Science and Technology Societies  (141S-1-3-0022).

\appendix
%
%
\section{Appendix} 
\setcounter{secnum}{\value {section}
\setcounter{equation}{0}
\renewcommand{\theequation}{\mbox{A.\arabic{equation}}}}

For $X=(x,t), Y=(y,s) \in \Bbb R^{n+1}$ we define the parabolic metric
$$d_p (X,Y)= \max\{ |x-y|, |s-t|^{\frac12} \},\quad X, Y \in \Bbb R^{n+1}.
$$
Let $Q=\Omega \times (a,b)$, where $\Omega \subset \Bbb R^n$ is a bounded domain, and $ -\infty <a<b<+\infty$. Then, for $0<\gamma <1$ we define the space of 
H\"{o}lder continuous functions on $Q$, 
$C^{\gamma, \frac{\gamma}{2}} (\bar{Q})$ 
by functions $f: \bar{Q} \to \Bbb R$ such that
$$
[f]_{C^{\gamma, \frac{\gamma}{2}}} =\sup_{X,Y\in \bar{Q}, X\not = Y}
 \frac{|f(X)-f(Y)|}{d_p (X,Y)^{\gamma} } <+\infty.
$$

\hspace{0.5cm}
The following parabolic  version of the Poincare inequality has been proved in 
\cite[Lemma\,B.3]{wol3}

\begin{lem}[Parabolic Poincar\'{e}-type inequality]
\label{lemA.1}
Let $Q_r=Q_r(X_0)\subset \R^{n+1}$\, $(n\in \N)$. 
Let $u \in  L^p(Q_r)$ be such that $\nabla u \in  L^p(Q_r)$\, $(1\le p<+\infty )$. 
In addition suppose that there exists
$\bbf \in L^1(Q_r)^n$ such that $\partial _t u= \nabla  \cdot \bbf$ in sense of distributions, i.\,e.
\be\label{A.1}
\intl_{Q_r} u \partial _t \varphi  dx dt = \intl_{Q_r} \bbf\cdot \nabla \varphi  dxdt   
\quad \forall\, \varphi \in C^\infty _{\rm c}(Q_r).
\ee
Then 
\be\label{A.2}
\intmw_{Q_r} |u- u_{Q_r} |^p dxdt \le c r^p \intmw_{Q_r} |\nabla u|^p dxdt 
+ c r^p\bigg(\intmw_{Q_r} |\bbf|dxdt\bigg)^p,      
\ee 
where $c=\const>0$, depending on $n$ and $p$ only, but not on $r, u$ or  $\bbf$.  
 
\end{lem}

\hspace{0.5cm}
The following elementary algebraic lemma has been used in  the proof of
 Theorem\,\ref{thm6.1}. 

\begin{lem}
\label{lemA.2}
Let $\{M_i\}$ and $ \{\lambda _i\}$ be  sequences of positive numbers such that 
$\suml \lambda _i\le  \frac {1} {2}$, and 
\be\label{A.10}
|M_{j+1}- M_{j}  |\le  (1+ M_j) \lambda _j\quad \forall\, j\in \N.
\ee
Then, 
\be\label{A.11}
M_i \le 1+ 2M_1\quad \forall\, i\in \N.
\ee
\end{lem}

{\it Proof } We prove the statement of this lemma by induction. Cleary, for $i=1$ the assertion is trivially fulfilled. Assume, \eqref{A.11} holds for $j=1,\ldots, i$. Then,  with help of of triangle inequality  
and \eqref{A.10} for  $j=1,\ldots,i$
we get 
\begin{align*}
M_{i+1} &\le M_1 + |M_{i+1} - M_1| \le M_1 + \suml_{j=1}^{i} |M_{j+1}- M_j| 
\\     
&\le M_1 + \suml_{j=1}^{i} (1+ M_j) \lambda _j \le M_1 +  (2 + 2 M_1) \suml_{j=1}^{i} \lambda _j
\le 1+ 2M_1.
\end{align*}
Whence, the claim is proved. \hfill \Beweisende


\begin{thebibliography}{99}

\bibitem{ach}
{\sc M.~Acheritogaray, P.~Degond, A.~Frouvelle, and J.-G. Liu}, {\em Kinetic
  formulation and global existence for the Hall-magnetohydrodynamic system},
  Kinetic and Related Models, 4 (2011), pp.~901--918.

\bibitem{caf}
{\sc L.~Caffarelli, R.~Kohn, and L.~Nirenberg}, {\em Partial regularity of
  suitable weak solutions of the Navier-Stokes equations}, Comm. Pure Appl.
  Math., 35 (1982), pp.~771--831.

\bibitem{cha1}
{\sc D.~Chae, P.~Degond, and J.-G. Liu}, {\em Well-posedness for
  Hall-magnetohydrodynamics}, Ann. Inst. Henri Poincare-Analyse Nonlineaire, 31
  (2014), pp.~555--565.

\bibitem{cha2}
{\sc D.~Chae and J.~Lee}, {\em On the blow-up criterion and small data global
  existence for the Hall-magnetohydrodynamics}, J. Differential Equations, 256
  (2014), pp.~3835--3858.

\bibitem{cha3}
{\sc D.~Chae and M.~Schonbek}, {\em On the temporal decay for the
  Hall-magnetohydrodynamic equations}, J. Differential Equations, 255 (2013),
  pp.~3971--3982.

\bibitem{cha4}
{\sc D.~Chae and J.~Wolf}, {\em On partial regularity for the steady 
Hall-magnetohydrodynamics system}, preprint,  (2015).

\bibitem{dum}
{\sc E.~Dumas and F.~Sueur}, {\em On the weak solutions to the
  Maxwell-Landau-Lifshitz equations and to the Hall-magnetohydrodynamic
  equations}, Comm. Math. Phys., 330 (2014), pp.~1179--1225.

\bibitem{fan}
{\sc J.~Fan, S.~Huang, and G.~Nakamura.}, {\em Well-posedness for the
  axisymmetric incompressible viscous Hall-magnetohydrodynamic equations},
  Appl. Math. Lett., 26 (2013), pp.~963--967.

\bibitem{for}
{\sc T.~Forbes}, {\em Magnetic reconnection in solar flares}, Geophys.
  Astrophys. Fluid Dyn., 62 (1991), pp.~15--36.

\bibitem{gia}
{\sc M.~Giaquinta}, {\em Multiple integrals in the calculus of variations and
  nonlinear elliptic systems}, ~Ann. of Math. Studies, vol. {\bf 105}, Princeton
  Univ. press, Princeton, New Jersey, 1983.

\bibitem{hom}
{\sc H.~Homann and R.~Grauer}, {\em Bifurcation analysis of magnetic
  reconnection in Hall-MHD systems}, Physica D, 208 (2005), pp.~59--72.

\bibitem{lad}
{\sc O.~A. Ladyzehnskaya and G.~A. Seregin}, {\em On partial regularity of
  sutiable weak solutions to the three-dimensional Navier-Stokes equations}, J.
  Math. Fluid Mech., 1 (1999), pp.~356--387.

\bibitem{lig}
{\sc M.~J. Lighthill}, {\em Studies on magnetohydrodynamic waves and other
  anisotropic wave motions}, Philos. Trans. R. Soc. Lond., Ser. A (1960),
  pp.~397--430.

\bibitem{lin}
{\sc F.~H. Lin}, {\em A new proof of the Caffarelli-Kohn-Nirenberg theorem},
  Comm. Pure Appl. Math., 51 (1998), pp.~241--257.

\bibitem{miu}
{\sc H.~Miura and D.~Hori}, {\em Hall effects on local structure in decaying
  MHD turbulence}, J. Plasma Fusion Res., 8 (2009), pp.~73--76.

\bibitem{pol}
{\sc J.~M. Polygiannakis and X.~Moussas}, {\em A review of magneto-vorticity
  induction in Hall- MHD plasmas}, Plasma Phys. Control \& Fusion, 43 (2001),
  pp.~195--221.

\bibitem{dap}
{\sc G.~D. Prato}, {\em Spazi ${\mathcal L}^{p, \theta}(\Omega, \delta)$ e loro
  propriet\`a}, Ann. Mat. Pura Appl., 69 (1965), pp.~383--392.

\bibitem{sch}
{\sc V.~Scheffer}, {\em Partial regularity of solutions to the Navier-Stokes
  equations}, Pacific J. Math., 66 (1976), pp.~535--552.

\bibitem{sha}
{\sc D.~Shalybkov and V.~Urpin}, {\em The Hall effect and the decay of magnetic
  fields}, Astron. Astrophys.,  (1997), pp.~685--690.

\bibitem{sim}
{\sc A.~N. Simakov and L.~Chac\'{o}n}, {\em Quantitative, analytical model for
  magnetic reconnection in Hall-magnetohydrodynamics}, Phys. Rev. Lett., 101
  (2008).

\bibitem{war}
{\sc M.~Wardle}, {\em Star formation and the Hall effect}, Astrophys. Space
  Sci., 292 (2004), pp.~317--323.

\bibitem{wol3}
{\sc J.~Wolf}, {\em "Regularit\"at schwacher L\"osungen elliptischer und
  parabolischer Systeme partieller Differentialgleichungen mit Entartung. Der
  fall $1 < p < 2$"}, Dissertation, Humboldt-Universit\"at zu Berlin, Berlin (2001).

\bibitem{wol1}
{\sc J.~Wolf}, {\em On the local regularity of suitable weak solutions to the
  generalized {Navier-Stokes} equations}, Annali della Universita Ferrara, (Doi
  10.1007/s11565-014-0203-6) (2014).

\bibitem{wol2}
{\sc J.~Wolf}, {\em On the local pressure of the Navier-Stokes equations and
  related systems}, submitted,  (2015).

\end{thebibliography}
\end{document}